\documentclass[12pt]{amsart}
\usepackage{epsfig}
\usepackage{graphics}
\numberwithin{equation}{section}
\theoremstyle{plain}
   \newtheorem{thm}{Theorem}[section]
   \newtheorem{cor}[thm]{Corollary}
   \newtheorem{prop}[thm]{Proposition}
   \newtheorem{lemma}[thm]{Lemma}

   \newtheorem{conj}[thm]{Conjecture}

\theoremstyle{definition}
   \newtheorem{dfn}[thm]{Definition}
   \newtheorem{ex}[thm]{Example}
\theoremstyle{remark}
   \newtheorem{rem}[thm]{Remark}

\newcommand{\m}{\underline{m}}

\newcommand{\N}{\mathbb N}

\newcommand{\R}{\mathbb R}

\newcommand{\Z}{\mathbb Z}

\newcommand{\ve}{\varepsilon}
\newcommand{\bd}{\partial}
\newcommand{\fM}{\mathfrak m}
\newcommand{\fL}{\mathfrak l}

\begin{document}
\title[Compatible contact structures of fibered Seifert links]{Compatible contact structures of fibered Seifert links in homology 3-spheres}
\author{Masaharu Ishikawa}
\footnote[0]{This work is partly supported by the Grant-in-Aid for Young Scientists (B), the Ministry of Education, Culture, Sports, Science and Technology, Japan.}
\address{Mathematical Institute\endgraf Tohoku University\endgraf Sendai 980-8578\endgraf Japan}
\email{ishikawa@math.tohoku.ac.jp}
\subjclass[2000]{Primary 57M50; Secondary: 55R25, 57R17, 32S55.}

\begin{abstract}
We study compatible contact structures of fibered Seifert multilinks
in homology $3$-spheres and especially give
a necessary and sufficient condition for the contact structure to be tight
in the case where the Seifert fibration is positively twisted.
As a corollary we determine the strongly quasipositivity of
fibered Seifert links in $S^3$.
We also study the compatible contact structures of cablings
along links in any $3$-manifolds.
\end{abstract}

\maketitle

\section{Introduction}

A contact structure on an closed, oriented, smooth $3$-manifold $M$
is the kernel of a $1$-form $\alpha$ on $M$ satisfying 
$\alpha\land d\alpha\ne 0$ everywhere. In this paper,
we only consider a positive contact form, i.e., a contact form 
$\alpha$ with $\alpha\land d\alpha>0$.
In~\cite{tw}, Thurston and Winkelnkemper used open book decompositions
to show the existence of contact structures on any $3$-manifolds.
In~\cite{giroux}, Giroux then focused on their idea, introduced 
the notion of contact structures supported by open book decompositions,
and studied the correspondence between contact structures up to contactomorphisms
and open book decompositions up to plumbings of positive Hopf bands, cf.~\cite{etnyre}.
Instead of the terminology ``supported'', we will say that
the contact structure is ``compatible'' with an open book decomposition and vice versa.

In the study of open book decompositions of $3$-manifolds,
it is important to determine if the compatible contact structure is tight or overtwisted
since it gives a rough classification of open book decompositions by Giroux's correspondence.
An explicit construction sometimes helps to determine the tightness.
For example, in~\cite{eo, ozbagci}, Etg\"{u} and Ozbagci gave explicit descriptions
of contact structures transverse to the fibers of circle bundles and certain
Seifert fibered manifolds and proved that such contact structures are Stein fillable.
Stein fillable contact structures are known to be tight by Eliashberg and Gromov~\cite{eliashberg:90, gromov}.

The purpose of this paper is to give an explicit construction of contact structures 
compatible with fibered Seifert links in homology $3$-spheres.
We hereafter use the terminology ``fibered link'' instead of ``open book decomposition''.
Following the book of Eisenbud and Neumann~\cite{en}, we denote a Seifert
fibered homology $3$-sphere as $\Sigma(a_1,a_2,\ldots,a_k)$,
where $a_i$'s are the denominators of the Seifert invariants.
The Seifert fibration has different properties depending on
the sign of the product $a_1a_2\cdots a_k$;
if $a_1a_2\cdots a_k>0$ then the fibers of the Seifert fibration
are twisted positively, as those of the positive Hopf fibration,
and if $a_1a_2\cdots a_k<0$ then they are negatively twisted.

A {\it Seifert link} $L$ in $\Sigma(a_1,\ldots,a_k)$ is an oriented link whose
exterior admits a Seifert fibration. 
Every Seifert link is realized as a union of a finite number 
of fibers of the Seifert fibration.
A multilink is a link each of whose link components is equipped with a non-zero integer,
called the {\it multiplicity}. 
A multilink is said to be {\it fibered} if its complement admits a fibration over $S^1$
such that the number of local leaves of 
the fiber surface in a small tubular neighborhood of each link component 
is the absolute value of its multiplicity and
the orientation induced from the fiber surface agrees with
the sign of the multiplicity, see Section~2 for precise definitions.
Note that a multilink is a usual link if all the multiplicities are in $\{-1,\; 1\}$. 
The criterion in~\cite[Theorem~11.1]{en} determines the fiberedness of 
a Seifert multilink in $\Sigma(a_1,\ldots,a_k)$, from which we can see that
Seifert multilinks are fibered in most cases.

Now we assign an orientation to the fibers of the Seifert fibration
under the assumption $a_1a_2\cdots a_k\ne 0$,
which we call the {\it orientation of the Seifert fibration}.
If the orientations of all the components of $L$ coincide with, or are opposite to, 
the orientation of the Seifert fibration
then we say that the orientation of $L$ is {\it canonical}.

In this paper we prove the following results.

\begin{thm}\label{thm01}
Let $L$ be a fibered Seifert multilink in $\Sigma(a_1,a_2,\ldots,a_k)$ with $a_1\cdots a_k>0$.
If the orientation of $L$ is canonical then the compatible contact structure is Stein fillable.
Otherwise it is overtwisted.
\end{thm}

The case $a_1a_2\cdots a_k<0$ will also be discussed in this paper.
As a consequence of our constructions in both cases,
we determine the tightness of fibered Seifert links in $S^3$.

\begin{figure}[htbp]
   \centerline{\input{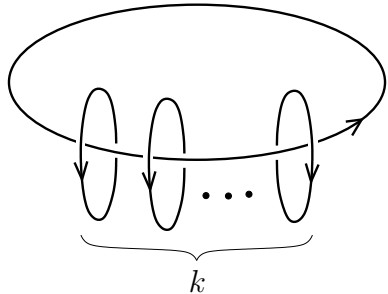}}
   \caption{Fibered Seifert links in case~(2).\label{fig111}}
\end{figure}

\begin{thm}\label{cor03}
Let $L$ be a fibered Seifert link in $S^3=\Sigma(a_1,a_2)$. 
Then the compatible contact structure of $L$ is tight
if and only if $L$ is one of the following cases:
\begin{itemize}
\item[(1)] $a_1a_2>0$ and the orientation of $L$ is canonical.
\item[(2)] $L$ is an oriented link described in Figure~\ref{fig111} with $k\geq 1$.
\end{itemize}
\end{thm}

With a small additional effort, we can remove 
the fiberedness assumption by replacing `tightness' into `strongly quasipositivity',
see Section~7 for the definition of strongly quasipositive links.

\begin{cor}\label{cor04}
Let $L$ be a non-splittable Seifert link in $S^3$.
Then, $L$ is strongly quasipositive if and only if it is in case~{\rm (1)} or {\rm (2)} above,
or in case~{\rm (3)} stated below:
\begin{itemize}
\item[(3)] $L$ is a negative torus link consisting of even number of link components
half of which have reversed orientation.
\end{itemize}
\end{cor}

Here a link $L$ in $S^3$ is called {\it splittable} if
$S^3\setminus L$ contains an incompressible $2$-sphere.
The only splittable Seifert links are trivial links with several components.

The technique of cabling with contact structure
can be used for studying cablings along fibered links in arbitrary $3$-manifolds.
Let $L(\m)$ be a fibered multilink in an oriented, closed, smooth $3$-manifold $M$ 
with cabling in a solid torus $N$ in $M$ and $L'(\m')$ be a fibered multilink 
obtained from $L(\m)$ by retracting $N$ into its core curve.
Note that $L'(\m')$ is always fibered.
We say that a cabling is {\it positive} if $L(\m)\cap N$ 
intersects the fiber surface of $L'(\m')$  positively transversely, 
and otherwise it is called {\it negative}.

\begin{thm}\label{thm04}
Let $L(\m)$ be a fibered multilink in an oriented, closed, smooth $3$-manifold $M$ 
with cabling in a solid torus $N$ in $M$ and $L'(\m')$ be the fibered multilink 
obtained from $L(\m)$ by retracting $N$ into its core curve.
Let $\xi$ and $\xi'$ denote the contact structures on $M$ compatible with 
$L(\m)$ and $L'(\m')$ respectively.
\begin{itemize}
\item[(1)] If $\xi'$ is tight and the cabling is positive, then $\xi$ is tight.
\item[(2)] If $\xi'$ is tight, the cabling is negative and $L(\m)\cap N$ has at least two components,
then $\xi$ is overtwisted.
\item[(3)] If $\xi'$ is tight, the cabling is negative, $L(\m)\cap N$ is connected, $p\geq 2$ and $q\leq -2$,
then $\xi$ is overtwisted.
\item[(4)] If $\xi'$ is overtwisted then $\xi$ is also overtwisted.
\end{itemize}
Here $p$ and $q$ are the coefficients of the slope $q\fM+p\fL$  of the cabling
with respect to the meridian-longitude pair $(\fM,\fL)$
on $\bd N$ which will be fixed in Section~8.1.
\end{thm}

The compatible contact structures of cablings in terms of multilinks are studied
independently by Baker, Etnyre and van Horn-Morris~\cite{behm}.
In their paper, a fibered multilink is called a rational open book decomposition.
The case of $M=S^3$ had been studied by Hedden in~\cite{hedden2}
using a different method.

This paper is organized as follows.
In Section~2, we fix the notation of Seifert fibered homology $3$-spheres 
and Seifert multilinks following the book~\cite{en}.
We introduce the notion of compatible contact structures for multilinks in Section~3.
The case $a_1\cdots a_k>0$ is studied in Section~4, including
the proof of Theorem~\ref{thm01},
and the case $a_1\cdots a_k<0$ is in Section~5, where
we give an explicit construction of contact structures and some criterion
for detecting overtwisted disks.
We then prove Theorem~\ref{cor03} in Section~6 and
Corollary~\ref{cor04} in Section~7.
In Section~8, we give the definitions of positive and negative cablings
and the proof of Theorem~\ref{thm04}.
A conjecture about strongly quasipositive orientation is posed
in the end of Section~7.

The author would like to thank Ko Honda, Shigeaki Miyoshi,
Jos\'{e} Mar\'{i}a Montesinos-Amilibia, Atsuhide Mori and
Kimihiko Motegi for their precious comments.

\section{Preliminaries}

In the following, $\text{int} X$ and $\bd X$ represent the interior and the boundary of a 
topological space $X$ respectively.

\subsection{Notation of Seifert fibered  homology $3$-spheres}

Let $\Sigma$ be a homology $3$-sphere.
We use the topological description of Seifert links in~\cite[p.60]{en}.
Let $\mathcal S=S^2\setminus\text{int}(D_1^2\cup\cdots\cup D_k^2)$
be a 2-sphere with $k$ holes and make
an oriented, closed, smooth $3$-manifold $\Sigma$ from $\mathcal S\times S^1$ by 
gluing solid tori $(D^2\times S^1)_1,\ldots,(D^2\times S^1)_k$
along the boundary $\bd (\mathcal S\times S^1)$.
To fix the notation, we first choose a section 
$\mathcal S^{\text{sec}}$ of $\pi:\mathcal S\times S^1\to \mathcal S$ and set
\[
\begin{split}
   Q_i&=(-\bd \mathcal S^{\text{sec}})\cap (D^2\times S^1)_i \\
   H&=\text{typical oriented fiber of $\pi$ in $\bd (D^2\times S^1)_i$}.
\end{split}
\]
Suppose that the gluing map of $(D^2\times S^1)_i$ to $\mathcal S\times S^1$ is 
given so that $a_iQ_i+b_iH$ is null-homologous in $(D^2\times S^1)_i$,
where $(a_i,b_i)\in\Z^2\setminus\{(0,0)\}$ and $\gcd(|a_i|,|b_i|)=1$.
To make the obtained $3$-manifold $\Sigma$ to be a homology $3$-sphere,
the integers $a_i$'s and $b_i$'s should satisfy the equality
$\sum_{i=1}^kb_ia_1\cdots a_{i-1} a_{i+1}\cdots a_k=\pm 1$.
Following~\cite{en}, in this paper, we always choose 
the coefficients $a_i$'s and $b_i$'s so that
$\sum_{i=1}^kb_ia_1\cdots a_{i-1} a_{i+1}\cdots a_k=1$
by replacing $(a_i,b_i)$ into $(-a_i,-b_i)$ for some $i$ if necessary.
Note that this equality ensures that
if one of $a_i$'s is zero then all the other $a_i$'s satisfy $|a_i|=1$,
and if $a_i\ne 0$ for all $i=1,\ldots,k$ then each pair $(i,j)$ with $i\ne j$ satisfies $\gcd(|a_i|,|a_j|)=1$.
Since the $3$-manifold $\Sigma$ does not depend on the ambiguity of the choice of $b_i$'s,
we may denote it as $\Sigma=\Sigma(a_1,\ldots,a_k)$.

The core curve $S_i$ of each solid torus $(D^2\times S^1)_i$
is a fiber of the Seifert fibration after the gluings.
We assign to $S_i$ an orientation
in such a way that the linking number of $S_i$ and $a_i Q_i+b_i H$ equals $1$.
This orientation is called the {\it working orientation}.

Let $(\fM_i,\fL_i)$ be the preferred meridian-longitude pair of the link complement 
$\Sigma\setminus S_i$ chosen such that the orientation of the longitude $\fL_i$ agrees
with the working orientation of $S_i$.
Then $(\fM_i,\fL_i)$ and $(Q_i,H)$ are related by the following equations, 
see~\cite[Lemma~7.5]{en}:
\begin{equation}\label{eqab}
   \begin{pmatrix} \fM_i \\ \fL_i \end{pmatrix}
   =
   \begin{pmatrix} a_i & b_i \\ -\sigma_i & \delta_i \end{pmatrix}
   \begin{pmatrix} Q_i \\ H \end{pmatrix}\quad\text{and}\quad
   \begin{pmatrix} Q_i \\ H \end{pmatrix}
   =
   \begin{pmatrix} \delta_i & -b_i \\ \sigma_i & a_i \end{pmatrix}
   \begin{pmatrix} \fM_i \\ \fL_i \end{pmatrix},
\end{equation}
where $\sigma_i=a_1\cdots\hat a_i\cdots a_k$ and 
$\delta_i=\sum_{j\ne i}b_ja_1\cdots \hat a_i\cdots \hat a_j\cdots a_k$. 
Note that they satisfy $a_i\delta_i+b_i\sigma_i=1$.

Set $A=a_1\cdots a_k$.
For a moment, we assume that $a_i\ne 0$ for all $i=1,\ldots,k$, in which case 
the orientation of the Seifert fibration in $\mathcal S\times S^1\to \mathcal S$
canonically extends into the fibers in $(D^2\times S^1)_i$ for each $i=1,\ldots,k$,
namely the orientation of the Seifert fibration of $\Sigma(a_1,\ldots,a_k)$
becomes well-defined.
Note that the working orientation on $S_i$ coincides with the orientation of
the Seifert fibration if and only if $a_i>0$.

\subsection{Fibered multilinks}

We give the definition of fibered multilinks in $3$-manifolds.
The same notion appears in~\cite{behm}, where 
the fibration is called a {\it rational open book decomposition}.

Let $M$ be an oriented, closed, smooth $3$-manifold
and $L$ an unoriented link in $M$ with $n$ link components.
We first assign an orientation to each link component of $L$,
which we also call a {\it working orientation}.
A {\it multilink} $L(\m)$ in $M$ is a link
each of whose components is equipped with a non-zero integer,
called the {\it multiplicity}, where $\m=(m_1,\ldots,m_n)$ represents
the set of multiplicities.
A multilink $L(\m)$ is called {\it fibered}
if there is a fibration $M\setminus L\to S^1$ such that 
\begin{itemize}
\item the intersection
of the fiber surface and a small tubular neighborhood $N(S_i)$ of each 
link component $S_i$ of $L(\m)$ locally 
consists of $|m_i|>0$ leaves meeting along $S_i$, and
\item the working orientation of $S_i$ is consistent 
with (resp. opposite to) the orientation induced from the fiber surface if $m_i>0$ (resp. $m_i<0$)
\end{itemize}
(cf.~\cite[p.28--29]{en}).

\subsection{Fibered Seifert multilinks}

A Seifert link $L$ in $\Sigma(a_1,\ldots,a_k)$ 
is a union of finite number of fibers of the Seifert fibration.
We had introduced the working orientation for each link component $S_i$ of $L$ in Section~2.1.
Using this working orientation, we assign a multiplicity to each $S_i$ 
and make $L$ to be a Seifert multilink. We denote this multilink as
\[
   L(\m)=(\Sigma(a_1,\ldots,a_k), m_1S_1\cup\cdots\cup m_nS_n),
\]
where $1\leq n\leq k$.
Note that Seifert multilinks are fibered in most cases and the fiberedness 
can be determined by a certain criterion stated in~\cite[Theorem~11.2]{en}. 
A typical example of non-fibered Seifert multilink is
the link obtained as the boundary of an $N$-times full-twisted annulus with $|N|\geq 2$.

Suppose that $L(\m)$ is fibered. 
The interiors of the fiber surfaces of $L(\m)$ intersect
the fibers of the Seifert fibration transversely
except for the case where $L(\m)$ is a positive or negative Hopf multilink,
see~\cite[Theorem~11.2]{en} and the proof therein.
In these exceptional cases, the transversality does not hold if the multiplicities
and the denominators of the Seifert invariants satisfy a certain equation.
As mentioned in~\cite[Proposition~7.3]{en}, a Seifert multilink is invertible
and this involution changes $L(\m)$ into $L(-\m)$.
In particular, this reverses the sign of the intersection of 
the interiors of the fiber surfaces and the fibers of the Seifert fibration.
So, by choosing one of $L(\m)$ and $L(-\m)$ suitably,
we often assume in this paper that the intersection is positive. 
We name it the {\it positive transverse property} and write it \text{\rm (PTP)} for short.

Now we consider the case where $A=a_1\cdots a_k\ne 0$.
In this case, as we already mentioned, the orientation of
the Seifert fibration of $\Sigma(a_1,\ldots,a_k)$ becomes well-defined. 

\begin{dfn}
Suppose $A\ne 0$. A link component $m_iS_i$ of a fibered Seifert multilink $L(\m)$ 
with \text{\rm (PTP)} is called
{\it positive} (resp. {\it negative}) if its orientation is consistent with
(resp. opposite to) the orientation of the Seifert fibration.
If the orientations of the link components of $L(\m)$ are either
all positive or all negative
then we say that the orientation of $L(\m)$ is {\it canonical}.
\end{dfn}

\section{Fibered multilinks and contact structures}

\subsection{A Lutz tube}

We first introduce terminologies in $3$-dimensional contact topology briefly,
see for instance~\cite{os, geiges} for general references.

\begin{figure}[htbp]
   \centerline{\input{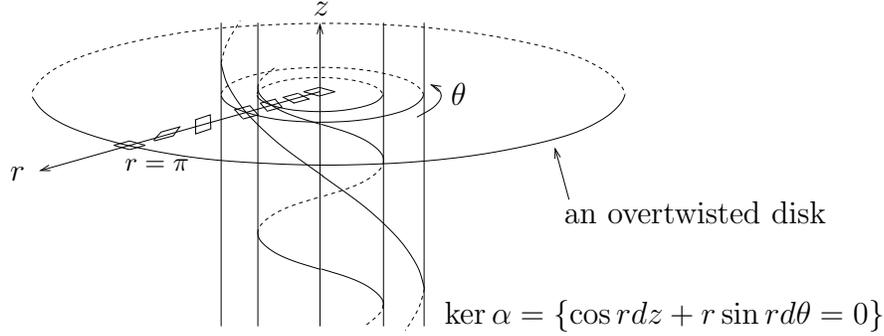}}
   \caption{A typical example of overtwisted contact structures.\label{fig16}}
\end{figure}

A {\it contact structure} on $M$ is the $2$-plane field given by
the kernel of a $1$-form $\alpha$ satisfying $\alpha\land d\alpha\ne 0$ everywhere on $M$.
In this paper, we only consider a contact structure given by
the kernel of a $1$-form $\alpha$ satisfying $\alpha\land d\alpha>0$,
called a {\it positive contact form} on $M$.
A vector field $R_\alpha$ on $M$ determined by the conditions 
$d\alpha(R_\alpha,\cdot)\equiv 0$ and $\alpha(R_\alpha)\equiv 1$ is called
the {\it Reeb vector field} of $\alpha$.
The $3$-manifold $M$ equipped with a contact structure $\xi$ is called
a {\it contact manifold} and denoted as $(M,\xi)$.
Two contact manifolds $(M_1,\xi_1)$ and $(M_2,\xi_2)$ are said to be
{\it contactomorphic} if there exists a diffeomorphism $\varphi:M_1\to M_2$
such that $d\varphi:TM_1\to TM_2$ satisfies $d\varphi(\xi_1)=\xi_2$.
A disk $D$ in $(M,\xi)$ is called {\it overtwisted} if 
$D$ is tangent to $\xi$ at each point on $\bd D$.
If $(M,\xi)$ has an overtwisted disk then we say that $\xi$ is {\it overtwisted}
and otherwise that $\xi$ is {\it tight}. 
A typical example of overtwisted contact structures is given as follows:
Let $\alpha$ be the contact form on $\R^3$ given by
\[
  \alpha=\cos rdz+r\sin rd\theta,
\]
where $(r,\theta,z)$ are coordinates of $\R^3$ with polar coordinates $(r,\theta)$.
The contact structure $\ker\alpha$ is as shown in Figure~\ref{fig16}.
We can find an overtwisted disk in the tube $\{(r,\theta,z) \,;\, |r|\leq \pi+\ve\}$,
where $\ve>0$ is a sufficiently small real number. 
Hence, this contact structure is overtwisted.

\begin{figure}[htbp]
   \centerline{\input{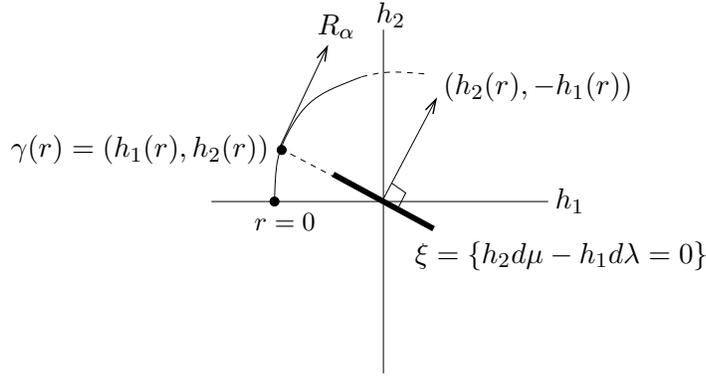}}
   \caption{How to read $\xi=\ker\alpha$ and $R_\alpha$ 
from the curve $\gamma(r)=(h_1(r),h_2(r))$.\label{fig1}}
\end{figure}

Now we introduce an effective way to describe a contact structure on $D^2\times S^1$. 
Let $\alpha$ be a $1$-form on $D^2\times S^1$ given by
$\alpha=h_2d\mu-h_1d\lambda$, where
$(r,\mu,\lambda)$ are coordinates of $D^2\times S^1$
with polar coordinates $(r,\mu)$ of $D^2$, and
$h_1$ and $h_2$ are real-valued smooth functions with parameter $r$.
We have
\[
\begin{split}
   d\alpha&=h'_2dr\land d\mu-h'_1dr\land d\lambda \\ 
   \alpha\land d\alpha&=(h_1'h_2-h_1h_2')dr\land d\mu\land d\lambda,
\end{split}
\]
where $h_1'$ and $h_2'$ are the derivatives of $h_1$ and $h_2$ with parameter $r$ respectively.
So, $\alpha$ is a positive contact form if and only if $h_1'h_2-h_1h_2'>0$.
We now plot $(h_1,h_2)$ on the $xy$-plane. Since $(h_2,-h_1)$ represents a vector normal to
the $2$-plane of the contact structure $\ker\alpha$, we can regard
the line connecting $(0,0)$ and $(h_1,h_2)$ as the slope of $\ker\alpha$.
The Reeb vector field $R_\alpha$ of $\alpha$ is given as
\[
   R_\alpha=\frac{1}{h_1'h_2-h_1h_2'}\left(h_1'\frac{\bd}{\bd \mu}+h_2'\frac{\bd}{\bd \lambda}\right).
\]
The parameter $r$ varies from $0$ to $1$, namely from $\{(0,0)\}\times S^1$ to the  boundary 
of $D^2\times S^1$, and the pair of functions $(h_1(r),h_2(r))$ defines a curve $\gamma$ on the $xy$-plane.
In summary, the curve $\gamma$ has the following properties:
\begin{itemize}
\item Since $h_1'h_2-h_1h_2'>0$, $(0,0)\not\in\gamma([0,1])$ and $\gamma$ moves in clockwise orientation.
\item The line connecting $(0,0)$ and $(h_1,h_2)$ represents the slope of $\ker\alpha$ and 
the vector $(h_2,-h_1)$ represents the positive side of $\ker\alpha$.
\item The speed vector $(h_1',h_2')$ is parallel to $R_\alpha$ and points in the same direction.
\end{itemize}
See Figure~\ref{fig1}. To make $\alpha$ to be
a well-defined contact form in a neighborhood of $r=0$,
we choose $\gamma$ near $r=0$ such that $(h_1,h_2)=(-c,r^2)$ or $(h_1,h_2)=(c,-r^2)$
with some positive constant $c$, so that 
$\alpha$ has the form $\alpha=r^2d\mu+cd\lambda$ or
$\alpha=-(r^2d\mu+cd\lambda)$ near $r=0$ respectively.

If the curve $\gamma$ intersects the positive $x$-axis, then the contact structure
$\ker\alpha$ on $D^2\times S^1$ has an overtwisted disk,
similar to Figure~\ref{fig16}, whose boundary corresponds to the intersection point 
of $\gamma$ and the positive $x$-axis.
In this paper, we call the tube $(D^2\times S^1, \ker\alpha)$ a {\it Lutz tube}
and use it frequently to show the existence of an overtwisted disk.

\subsection{Contact structures compatible with multilinks}

The notion of compatible contact structures of fibered links can be generalized
to fibered multilinks canonically. This idea also appears in~\cite{behm}.
Let $M$ be a closed, oriented, smooth $3$-manifold.

\begin{dfn}
A fibered multilink $L(\m)$ in $M$ is said to be {\it compatible} with 
a contact structure $\xi=\ker\alpha$ on $M$ if $L(\m)$ is positively transverse to $\xi$ and
$d\alpha$ is a volume form on the interiors of the fiber surfaces of $L(\m)$.
\end{dfn}

The next lemma gives a useful interpretation of the notion of compatible
contact structures
in terms of Reeb vector fields. In this paper we mainly use this characterization.

\begin{lemma}\label{lemma_etnyre}
A fibered multilink $L(\m)$ in $M$ is compatible with a contact structure $\xi$ on $M$ 
if and only if there exists a contact form $\alpha$ on $M$
with $\xi=\ker\alpha$ such that the Reeb vector field $R_\alpha$
is tangent to $L(\m)$ and positively transverse to the interiors of the fiber surfaces of $L(\m)$,
and its orientation is consistent with that of $L(\m)$ induced from the fiber surfaces.
\end{lemma}

\begin{proof}
The proof for a fibered link in~\cite[Lemma~3.5]{etnyre} works in this case also.
\end{proof}

Now we introduce two fundamental facts concerning compatible contact structures
of fibered multilinks, following the fibered link case.

\begin{prop}\label{lemmatw}
Any fibered multilink in $M$ admits a compatible contact structure.
\end{prop}

Although the proof is analogous to the one in~\cite{tw},
since an explicit contact form of the compatible contact structure 
will be needed in the proof of Lemma~\ref{thm3} later, we prove the assertion here
with presenting the contact form.
A similar proof can be found in~\cite{behm}.

\begin{proof}
Let $L(\m)$ be a fibered multilink in $M$ with
$n$ link components $m_1S_1,\ldots, m_nS_n$ and
$N(S_i)$ a small compact tubular neighborhood of $S_i$ in $M$ for $i=1,\ldots,n$.
We denote by $F_t$ the fiber surface of $L(\m)$ over $t\in S^1=[0,1]/0\sim 1$ and
choose a diffeomorphism $\phi_t:F_0\to F_t$ of the fibration of $L(\m)$
in such a way that
\[
\phi_t(r_i,\mu_i,\lambda_i)=\left(r_i, \mu_i+\frac{t}{|m_i|} , \lambda_i \right)
\]
in $N(S_i)$, where $(r_i,\mu_i,\lambda_i)$ are coordinates of $N(S_i)=D^2\times S^1$
chosen such that $(r_i,\mu_i)$ are the polar coordinates of $D^2$ and
the orientation of $\lambda$ agrees with that of
the corresponding link component of $L(\m)$.
For convenience, we set the coordinates $(r_i,\mu_i)$ such that the radius of $D^2$ is $1$.

Let $\theta_i$ be the coordinate function on the curve $-(F_0\setminus\text{int\,}N(S_i))$
given as $\theta_i=-\lambda_i$.
Then, as in~\cite{tw}, we can find a $1$-form $\beta$ on 
$F_0\cap (\mathcal S\times S^1)$ such that $d\beta$ is a volume form on $F_0\cap(\mathcal S\times S^1)$
and $\beta=-(1/r_i)d\theta_i$ near $\bd N(S_i)$.
The manifold $M$ is constructed from 
$F_0\times [0,1]$ by identifying $(x,1)\sim (\phi_1(x),0)$ for each $x\in F_0$
and then filling the boundary components by the solid tori $N(S_i)$'s.
According to this construction, we define a $1$-form $\alpha_0$ on $\mathcal S\times S^1$ as
\[
   \alpha_0=(1-t)\beta+t\phi_1^*(\beta)+Rdt,
\]
with $R>0$, which is given near $\bd N(S_i)$ as
\begin{equation}\label{eqsec3}
   \alpha_0=-\frac{1}{r_i}d\theta_i+Rdt=\frac{1}{r_i}d\lambda_i+R(v_id\mu_i-u_id\lambda_i),
\end{equation}
where $(u_i,v_i)$ is a vector representing
the oriented boundary of $F_0\setminus\text{int\,}N(S_i)$ 
on $\bd N(S_i)$ with coordinates $(\mu_i,\lambda_i)$; in other words,
$(v_i, -u_i)$ is a vector positively normal to $F_0$ on $\bd N(S_i)$.
Note that $v_i>0$. We choose $R$ sufficiently large so that $\alpha_0$ becomes 
a positive contact form on $\mathcal S\times S^1$.

For each $N(S_i)$, 
we extend $\alpha_0$ into $N(S_i)$ by describing a curve $\gamma(r_i)$
on the $xy$-plane explained in Section~3.1.
The endpoint $(h_1(1),h_2(1))$ of $\gamma(r_i)$ is given as $(h_1(1),h_2(1))=(Ru_i-1, Rv_i)$
and the speed vector $\gamma'(r_i)$ at $r_i=1$ is $(h_1'(r_i),h_2'(r_i))=(1,0)$. So,
we can describe a curve $\gamma(r_i)$ representing a positive contact form on $N(S_i)$ such that 
\begin{itemize}
\item $(h_1,h_2)=(-c, r^2)$ near $r=0$ with $c>0$,
\item $\gamma(1)$ and $\gamma'(1)$ satisfy the above conditions, and
\item $\gamma'(r_i)$ rotates monotonously.
\end{itemize}
Thus the contact form $\alpha_0$ is extended into $N(S_i)$.
We denote the obtained contact form on $M$ as $\alpha$.

Since the fibers of the Seifert fibration intersect $F_t\cap (\mathcal S\times S^1)$ 
positively transversely, $\ker\alpha$ is compatible with $L(\m)$ on $\mathcal S\times S^1$.
In each $N(S_i)$, we can isotope $F_t$ into the position shown in Figure~\ref{fig4c}
such that $\ker\alpha$ is compatible with $L(\m)$. This completes the proof.
\begin{figure}[htbp]
   \centerline{\input{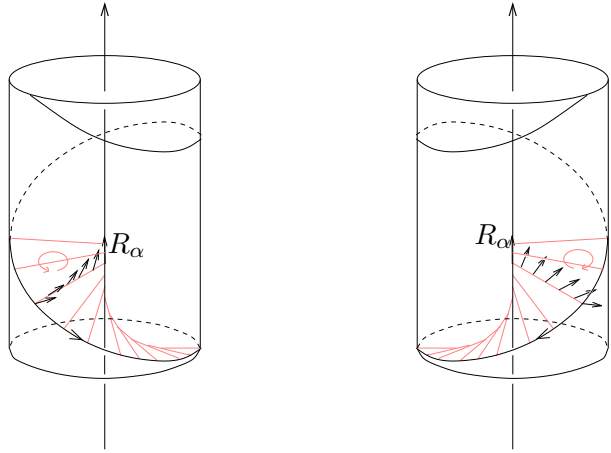}}
   \caption{The compatibility in the neighborhood $N(S_i)$.\label{fig4c}}
\end{figure}
\end{proof}

\begin{prop}\label{lemmatw2}
If two contact structures on $M$ are compatible with
the same fibered multilink in $M$ then they are contactomorphic.
\end{prop}

\begin{proof}
The proof for a fibered link in~\cite{giroux} works in this case also (cf.~\cite[Proposition~9.2.7]{os}).
\end{proof}

\section{Case $a_1a_2\cdots a_k>0$}

\subsection{Explicit construction of the contact structure}

Throughout this section, we always assume that $A=a_1\cdots a_k>0$.
Theorem~\ref{thm01} follows from the explicit construction of compatible
contact structures described below.

\begin{prop}\label{thm1}
Let $L(\m)=(\Sigma, m_1S_1\cup\cdots\cup m_nS_n)$
be a fibered Seifert multilink in a homology $3$-sphere $\Sigma=\Sigma(a_1,\ldots,a_k)$ 
with $A>0$. Assume \text{\rm (PTP)}.
Then there exists a positive contact form $\alpha$ on $\Sigma$
with the following properties:
\begin{itemize}
\item[(1)] $L(\m)$ is compatible with the contact structure $\xi=\ker\alpha$.
\item[(2)] The Reeb vector field $R_\alpha$ of $\alpha$ 
is tangent to the fibers of the Seifert fibration on $\mathcal S\times S^1$.
\item[(3)] The neighborhood $(D^2\times S^1)_i$ of each negative component $m_iS_i$ of $L(\m)$
 contains a Lutz tube. In particular, it contains an overtwisted disk.
\item[(4)] On the other $(D^2\times S^1)_i$'s, 
$\ker\alpha$ is transverse to the fibers of the Seifert fibration.
\end{itemize}
\end{prop}

\begin{rem}
The most canonical way to construct a contact structure compatible with a given fibered link
is to use the fiber surface as done in~\cite{tw}.
However this is difficult in our situation
because there is no systematic way to describe the fiber surface.
The idea of the proof of Theorem~\ref{thm01} is that we choose the contact form
such that its Reeb vector field is tangent to the fibers of the Seifert fibration everywhere except
in small neighborhoods of the negative components.
This makes sure that the contact structure is compatible with 
the fibered multilink in the most part.
The rest is done by describing possible local positions of the fiber surfaces
along the exceptional components.
\end{rem}

\begin{rem}
The existence of $S^1$-invariant contact forms on orientable
Seifert fibered $3$-manifolds is known in~\cite{kt}. 
The existence of a contact structure transverse to the fibers of a Seifert fibration
had been studied in~\cite{st} for circle bundles over closed surfaces and
in~\cite{lm} for Seifert fibered $3$-manifolds. 
The transverse contact structures are always Stein fillable as mentioned in~\cite[Theorem~4.2]{ch},
cf.~\cite{eo, ozbagci}.
This fact will be used in the proof of Theorem~\ref{thm01}.
\end{rem}

To prove Proposition~\ref{thm1}, we apply the argument in the proof in~\cite{tw} 
to the Seifert fibration. 
We denote the boundary component $(-\bd \mathcal S)\cap D^2_i$ of $\mathcal S$ by $C_i$.

\begin{lemma}\label{lemma1}
Suppose $A>0$ and let $U_i$ be a collar neighborhood of $C_i$ in $\mathcal S$ 
with coordinates $(r_i,\theta_i)\in [1,2)\times S^1$
satisfying $\{(r_i,\theta_i) \,;\, r_i=1\}=C_i$.
Then there exists a $1$-form $\beta$ on $\mathcal S$ which satisfies the following properties:
\begin{itemize}
\item[(1)] $d\beta>0$ on $\mathcal S$.
\item[(2)] If $b_i/a_i\leq 0$ then 
$\beta=R_ir_id\theta_i$ with $-b_i/a_i<R_i$ near $C_i$ on $U_i$.
\item[(3)] If $b_i/a_i>0$ then 
$\beta=(R_i/r_i)d\theta_i$ with $-b_i/a_i<R_i<0$ near $C_i$ on $U_i$.
\end{itemize}
\end{lemma}

\begin{proof}
Since $\sum_{i=1}^k(-b_i/a_i)=-1/A<0$,
we can choose $R_1,\ldots,R_k$ such that they satisfy the inequalities in~(2) and~(3) and 
the inequality $\sum_{i=1}^k R_i<0$.
Let $\Omega$ be a volume form on $\mathcal S$ which satisfies
\begin{itemize}
\item $\int_{\mathcal S}\Omega=-\sum_{i=1}^k R_i>0$,
\item $\Omega=R_idr_i\land d\theta_i$ near $C_i$ with $b_i/a_i\leq 0$, and
\item $\Omega=-(R_i/r_i^2)dr_i\land d\theta_i$ near $C_i$ with $b_i/a_i> 0$.
\end{itemize}
Let $\eta$ be any $1$-form on $\mathcal S$ which equals $R_ir_id\theta_i$ 
if $b_i/a_i\leq 0$ and $(R_i/r_i)d\theta_i$ if $b_i/a_i> 0$
near $C_i$.
By Stokes' theorem, we have
\[
\begin{split}
   \int_{\mathcal S}(\Omega-d\eta)&=\int_{\mathcal S}\Omega-\int_{\bd \mathcal S}\eta
=\int_{\mathcal S}\Omega+\sum_{i=1}^k\int_{C_i}R_id\theta_i \\
&=\int_{\mathcal S}\Omega+\sum_{i=1}^kR_i=0.
\end{split}
\]
Here $C_i$ is oriented as $-\bd \mathcal S$.
The closed $2$-form $\Omega-d\eta$ represents the trivial class in cohomology
vanishing near $\bd \mathcal S$. By de Rham's theorem, there is a $1$-form $\gamma$
on $\mathcal S$ vanishing near $\bd \mathcal S$ and satisfying $d\gamma=\Omega-d\eta$.
Define $\beta=\eta+\gamma$, then $d\beta=\Omega$ is a volume form on $\mathcal S$
and $\beta$ satisfies properties~(2) and~(3) near $\bd \mathcal S$ as required.
\end{proof}

We prepare two further lemmas which will be used for constructing the contact form
on $(D^2\times S^1)_i$.
Let $B=[1,2)\times S^1\times S^1\subset \mathcal S\times S^1$ be
a neighborhood of a boundary component of $\mathcal S\times S^1$ with coordinates $(r,\theta,t)$. 
We glue $D^2\times S^1$ to $B$ as
\[
   \mu \fM+\lambda \fL=(a\mu-\sigma\lambda)Q+(b\mu+\delta\lambda)H,
\]
where $(\fM,\fL)$ is a standard meridian-longitude pair of 
$\bd D^2\times S^1\subset D^2\times S^1$,
$Q$ is the oriented curve given by $\{1\}\times S^1\times\{\text{a point}\}\subset \bd B$,
$H$ is a typical fiber of the projection $[1,2)\times S^1\times S^1\to [2,1)\times S^1$ 
which omits the third entry, and $a, b, \sigma, \delta\in\Z$ are given 
according to relations~\eqref{eqab}.
The fibers $H=\sigma \fM+a\fL$ of the Seifert fibration on $\bd D^2\times S^1$
are canonically extended to the interior of $D^2\times S^1$.

\begin{lemma}\label{extension1}
Suppose $a\ne 0$ and either {\rm (i)}~$0\leq -b/a<R$ and $\alpha_0=Rrd\theta+dt$ or 
{\rm (ii)}~$-b/a<R<0$ and $\alpha_0=(R/r)d\theta+dt$, 
where $\alpha_0$ is a contact form on $B$.
Then there exists a contact form $\alpha$ on $B\cup (D^2\times S^1)$ 
with the following properties:
\begin{itemize}
\item[(1)] $\alpha=\alpha_0$ on $B$.
\item[(2)] $\ker\alpha$ is transverse to the fibers of the Seifert fibration in $D^2\times S^1$.
\item[(3)] $R_\alpha$ is tangent to $\{(0,0)\}\times S^1$
and the direction of $R_\alpha$ is consistent with the orientation of the Seifert fibration.
\item[(4)] $R_\alpha$ rotates monotonously with respect to the parameter $r\in [0,1]$.
\end{itemize}
\end{lemma}

\begin{proof}
We consider case~(i). 
Let $\sigma$ and $\delta$ be integers satisfying relations~\eqref{eqab}. 
Denote the gluing map of $D^2\times S^1$ to $B$ by $\varphi$, then we have
\begin{equation}\label{eqaaa}
\begin{split}
   \varphi^*\alpha_0&=Rrd(a\mu-\sigma\lambda)+d(b\mu+\delta\lambda)=(b+aRr)d\mu+(\delta-\sigma Rr)d\lambda \\
   &=a\left(\frac{b}{a}+Rr\right)d\mu+\frac{1}{a}\left(1-a\sigma\left(\frac{b}{a}+Rr\right)\right)d\lambda.
\end{split}
\end{equation}
If $a>0$ then $a(b/a+Rr)>0$ near $r=1$.
So, on the $xy$-plane, the point $(h_1(1),h_2(1))$ lies in the region $y>0$.
Since $R_{\alpha_0}$ is positively transverse to $\ker\alpha_0$ at $r=1$,
we can describe a smooth curve $\gamma(r)=(h_1(r),h_2(r))$ on the $xy$-plane 
representing a positive contact form on $B\cup(D^2\times S^1)$ such that
\begin{itemize}
\item $(h_1,h_2)=(-c,r^2)$  near $r=0$ with $c>0$, 
\item $h_2d\mu-h_1d\lambda=\varphi^*\alpha_0$ near $r=1$, and
\item $\gamma'(r)$ rotates monotonously,
\end{itemize}
as shown in Figure~\ref{fig2}. This satisfies the required properties.
\begin{figure}[htbp]
   \centerline{\input{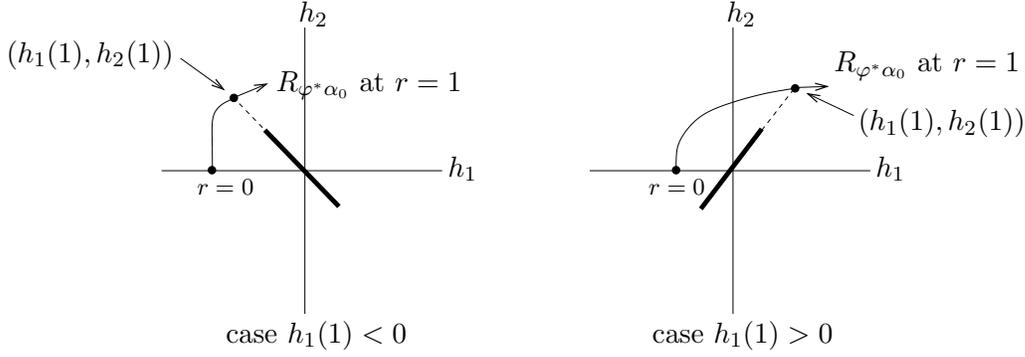}}
   \caption{Curves representing contact forms on $D^2\times S^1$ in 
Lemma~\ref{extension1}. The figures are in case $a>0$.\label{fig2}}
\end{figure}

If $a<0$ then $a(b/a+Rr)<0$ near $r=1$
and hence the point $(h_1(1),h_2(1))$ lies in the region $y<0$.
We choose a smooth curve $\gamma(r)$ such that
\begin{itemize}
\item $(h_1,h_2)=(c,-r^2)$ near $r=0$ with $c>0$, 
\item $h_2d\mu-h_1d\lambda=\varphi^*\alpha_0$ near $r=1$, and
\item $\gamma'(r)$ rotates monotonously.
\end{itemize}
Note that such a curve $\gamma(r)$ is given by the $\pi$-rotation of the figures
in Figure~\ref{fig2}.
The contact form $\alpha$ on $B\cup (D^2\times S^1)$ defined by this curve
satisfies the required properties as before.

The proof for case~(ii) is similar.
\end{proof}

\begin{lemma}\label{extension2}
Let $\alpha_0$ be a contact form on $B$ given by either {\rm (i)}~$\alpha_0=Rrd\theta+dt$ 
with $R>0$ or {\rm (ii)}~$\alpha_0=(R/r)d\theta+dt$ with $R<0$.
Then there exists a contact form $\alpha$ on $B\cup (D^2\times S^1)$ 
with the following properties:
\begin{itemize}
\item[(1)] $\alpha=\alpha_0$ on $B$.
\item[(2)] $\ker\alpha$ is transverse to the fibers of the Seifert fibration in $D^2\times S^1$
except on a torus $\{r_1\}\times S^1\times S^1$ embedded in $D^2\times S^1$ for some $r_1\in (0,1)$.
\item[(3)] $R_\alpha$ is tangent to $\{(0,0)\}\times S^1$ 
and the direction of $R_\alpha$ is opposite to the orientation of the Seifert fibration.
\item[(4)] $R_\alpha$ rotates monotonously with respect to the parameter $r\in [0,1]$.
\end{itemize}
Furthermore, if $R$ satisfies $R>-b/a$
then $(D^2\times S^1, \ker\alpha)$ contains a Lutz tube.
\end{lemma}

\begin{figure}[htbp]
   \centerline{\input{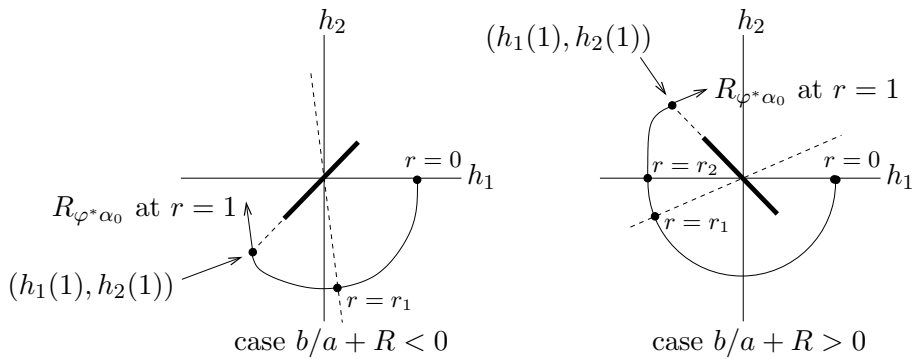}}
   \caption{Curves representing contact forms on $D^2\times S^1$
in Lemma~\ref{extension2}.\label{fig3}}
\end{figure}

\begin{proof}
The proof is analogous to the proof of Lemma~\ref{extension1}.
In case (i) with $a>0$, we choose a curve $\gamma$ on the $xy$-plane such that
 $(h_1,h_2)=(c,-r^2)$ near $r=0$ with $c>0$ as shown in Figure~\ref{fig3}.
This satisfies the required properties.
If $R>-b/a$ then a Lutz tube appears at $r=r_2$ as described
on the right in the figure.
The proofs in case $a<0$ and case (ii) are similar.
\end{proof}

\noindent
{\it Proof of Proposition~\ref{thm1}}.\;\;
Let $\alpha_0$ be the $1$-form on $\mathcal S\times S^1$ defined by $\alpha_0=\beta+dt$,
where $\beta$ is a $1$-form constructed
in Lemma~\ref{lemma1} and $t$ is the coordinate of $S^1$,
which is assumed to be consistent with the orientation of the Seifert fibration.
Since $\beta\land d\beta$ is a $3$-form on $\mathcal S$, we have
$\beta\land d\beta=0$ and
\[
   \alpha_0\land d\alpha_0=\beta\land d\beta+dt\land d\beta=d\beta\land dt>0.
\]
Thus $\alpha_0$ is a positive contact form on $\mathcal S\times S^1$
and its Reeb vector field is given by
$R_{\alpha_0}=\bd/\bd t$.
Note that, since $R_{\alpha_0}$ is tangent to the fibers of 
$\pi:\mathcal S\times S^1\to \mathcal S$ in the same direction, 
(PTP) implies that $R_{\alpha_0}$ is 
positively transverse to the fiber surfaces of $L(\m)$ in $\mathcal S\times S^1$.

Now we extend ${\alpha_0}$ into $(D^2\times S^1)_i$ in the following way.
If either $m_iS_i$ is a positive component or $i>n$ 
then we use the construction of a contact form in Lemma~\ref{extension1},
otherwise we use the construction in Lemma~\ref{extension2}.
We denote the extended contact form on $\Sigma$ by $\alpha$.

From the construction, we only need to check property~(1) in the assertion.
Due to Lemma~\ref{lemma_etnyre}, it is enough to check if $R_\alpha$ is tangent
to $L(\m)$ in the same direction and positively transverse to the interiors of 
the fiber surfaces of $L(\m)$.
This positive transversality had already been established in $\mathcal S\times S^1$.

\begin{figure}[htbp]
   \centerline{\input{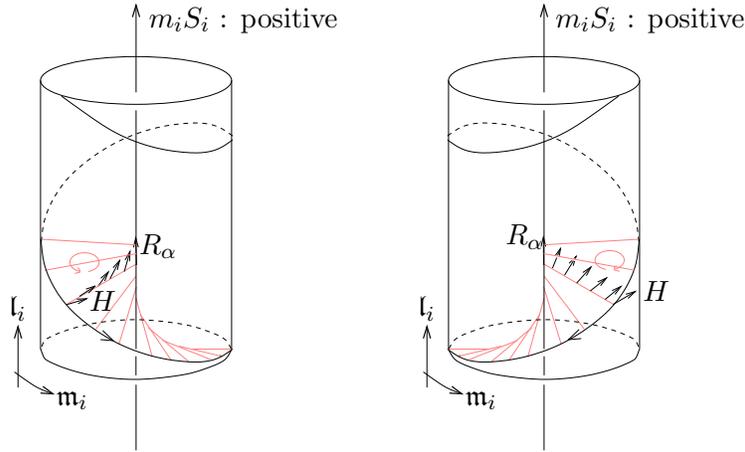}}
   \caption{The compatibility in the neighborhood $(D^2\times S^1)_i$ of a
positive component $m_iS_i$ in case $a_i>0$.\label{fig4a}}
\end{figure}

We first check the positive transversality in the neighborhood
$(D^2\times S^1)_i$ of a positive component $m_iS_i$.
Figure~\ref{fig4a} shows the mutual positions of the fiber surfaces $F$, the oriented 
fibers $H$ of the Seifert fibration and the Reeb vector field $R_\alpha$ on
$(D^2\times S^1)_i$ in case $a_i>0$.
The orientations of the link component $m_iS_i$ and the fibers $H$ are
as shown in the figures since $m_iS_i$ is a positive component, $a_i>0$, 
$\sigma_i>0$, and $H$ is given as $H=\sigma_i\fM_i+a_i\fL_i$.
The Reeb vector field $R_\alpha$ had already been given in the above construction.
Now there are three possibilities of the framing of the fiber surface $F$,
namely it is either positive, negative, or parallel to $m_iS_i$.
The case of positive framing is described on the left in the figure
and the case of negative framing is on the right. The parallel case is omitted.
In either case, we can isotope the fiber surfaces $F$ in $(D^2\times S^1)_i$
such that it satisfies the property~(1). Note that
the vectors of $R_\alpha$ on the right figure are directed under the fiber surface.
The proof in case $a_i<0$ is similar, in which case the figures are those
in Figure~\ref{fig4a} with replacing $(\fM_i, \fL_i)$ by $(-\fM_i, -\fL_i)$.

The property~(1) in $(D^2\times S^1)_i$ with $i>n$ can also be 
checked from the figure because the fiber surfaces on $(D^2\times S^1)_i$
consists of horizontal disks.

\begin{figure}[htbp]
   \centerline{\input{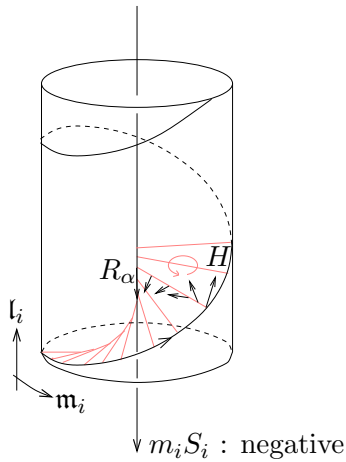}}
   \caption{The compatibility on the neighborhood $(D^2\times S^1)_i$ of a
negative component $m_iS_i$ in case $a_i>0$.\label{fig4b}}
\end{figure}

Suppose that $m_iS_i$ is a negative component. We assume that $a_i>0$.
Then the orientations of the link component $m_iS_i$ and the fibers $H$
become as shown in Figure~\ref{fig4b}.
There is only one possibility of the framing of the fiber surface $F$,
which is shown in the figure, otherwise they do not satisfy (PTP) on the boundary of
$(D^2\times S^1)_i$.
As shown in the figure, we can isotope the fiber surface $F$ in $(D^2\times S^1)_i$
such that it satisfies the property~(1).
The proof in case $a_i<0$ is similar and
the figure is as in Figure~\ref{fig4b} with replacing $(\fM_i, \fL_i)$ by $(-\fM_i, -\fL_i)$.
\qed
\vspace{3mm}

\subsection{Proof of Theorem~\ref{thm01}}

The next lemma will be used in the proof of Theorem~\ref{thm01}.

\begin{lemma}\label{lemma19}
If $A>0$ then every fibered Seifert multilink has at least one positive component.
\end{lemma}

\begin{proof}
Let $F$ be a fiber surface of a fibered Seifert multilink $L(\m)$ and
assume that $L(\m)$ has no positive component.
The fibers of the Seifert fibration are given as $H=\sigma_i \fM_i+a_i\fL_i$, 
where $\sigma_ia_i=A>0$.
Let $\gamma_i=u_i\fM_i+v_i\fL_i$ be the oriented boundary 
$\bd (F\cap (D^2\times S^1)_i)\setminus m_iS_i$,
where $u_i\in\Z$ and $v_i\in\Z\setminus\{0\}$ are chosen such that
the number of connected components of $\bd (F\cap (D^2\times S^1)_i)\setminus m_iS_i$
is equal to $\gcd(|u_i|,|v_i|)$ in case $u_i\ne 0$ and $|v_i|$ otherwise.
From (PTP), we have the inequality $I(\gamma_i,H)=u_ia_i-v_i\sigma_i>0$,
where $I(\gamma_i,H)$ is the algebraic intersection number of $\gamma_i$ and $H$ on $\bd (D^2\times S^1)_i$.
Furthermore, the fiber surface $F$ along $m_iS_i$ is given as shown in Figure~\ref{fig5} and
we can verify the inequality $a_iv_i>0$ from these figures.

\begin{figure}[htbp]
   \centerline{\input{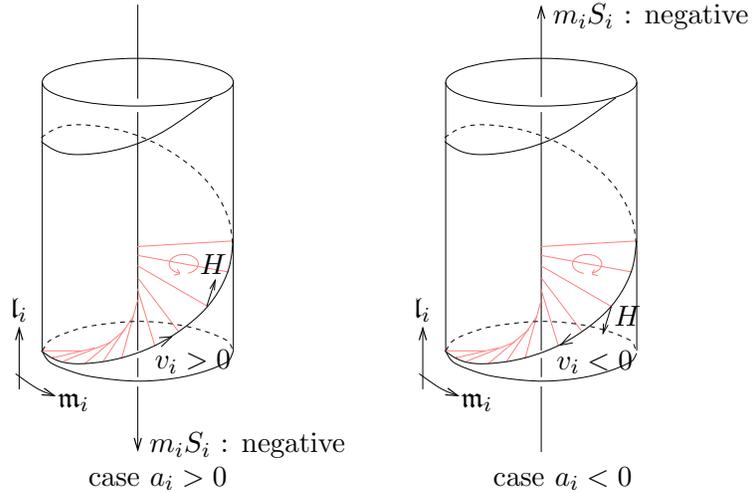}}
   \caption{The framing of $F$ along $m_iS_i$.\label{fig5}}
\end{figure}

For each $i=1,\cdots,n$, 
\[
   u_i\fM_i+v_i\fL_i=(a_iu_i-\sigma_iv_i)Q_i+(b_iu_i+\delta_iv_i)H.
\]
The union of these curves is homologous to the boundary of the fiber surface 
because it is a Seifert surface,
and hence the sum $\sum_{i=1}^n(u_i\fM_i+v_i\fL_i)$ is null-homologous in 
the complement $\Sigma\setminus L(\m)$.
This complement is obtained from $\mathcal S\times S^1$ by 
gluing $(D^2\times S^1)_i$, for $i=n+1,\ldots,k$,
in such a way that $a_iQ_i+b_iH$ corresponds to the meridian of $(D^2\times S^1)_i$. Hence
there exists a non-zero vector $(w_{n+1},\cdots,w_k)$ which satisfies
\[
   \sum_{i=1}^n ((a_iu_i-\sigma_iv_i)Q_i+(b_iu_i+\delta_iv_i)H)+\sum_{i=n+1}^k w_i(a_iQ_i+b_iH)=0.
\]
Since $\sum_{i=1}^k Q_i=0$ in $H_1(\mathcal S\times S^1)$ is the unique relation which we can use for
vanishing the coefficients of $Q_i$'s, all coefficients of $Q_i$'s must be the same value.
Hence we have the equality
\[
   \sum_{i=1}^n \left(Q_i+\frac{b_iu_i+\delta_iv_i}{a_iu_i-\sigma_iv_i}H\right)
+\sum_{i=n+1}^k \left(Q_i+\frac{b_i}{a_i}H\right)=0,
\]
which implies
\begin{equation}\label{eq100}
\begin{split}
   0&=\sum_{i=1}^n \frac{b_iu_i+\delta_iv_i}{a_iu_i-\sigma_iv_i}+\sum_{i=n+1}^k \frac{b_i}{a_i} 
     =\sum_{i=1}^n \left(\frac{b_i}{a_i}+\frac{v_i}{a_i(a_iu_i-\sigma_iv_i)}\right)
     +\sum_{i=n+1}^k \frac{b_i}{a_i} \\
    &=\frac{1}{A}+\sum_{i=1}^n\frac{v_i}{a_i(a_iu_i-\sigma_iv_i)}.
\end{split}   
\end{equation}
However the right hand side of this equation must be strictly positive since 
$a_iu_i-\sigma_iv_i>0$ and $a_iv_i>0$, which is  a contradiction.
\end{proof}

\noindent
{\it Proof of Theorem}~\ref{thm01}.
We first remark that it is enough to observe the tightness for
a specific contact form whose contact structure is compatible with $L(\m)$ by Proposition~\ref{lemmatw2}.
Assume that $L(\m)$ is not a Hopf multilink in $S^3$.
If all components of $L(\m)$ are negative
then it does not satisfy (PTP) by Lemma~\ref{lemma19}.
So, in this case, we reverse the orientation of $L(\m)$ as $L(-\m)$ so that
all components become positive.
If all components of $L(\m)$ are positive, then the compatible contact structure constructed 
according to the recipe in Proposition~\ref{thm1} is positively transverse
to the fibers of the Seifert fibration everywhere.
In particular, it is known that such a contact structure is always tight,
see~\cite{mw} and~\cite[Corollary~2.2]{lm}.
Moreover, since the monodromy of the fibration of $L(\m)$ is periodic, 
we can conclude that the contact structure is Stein fillable,
see~\cite[Theorem~4.2]{ch}.

Suppose that $L(\m)$ has at least one positive component and one negative component.
In this case, even if we reverse the orientation of $L(\m)$ by involution,
$L(\m)$ still has a negative component. Therefore, in either case,
the contact structure $\ker\alpha$ has an overtwisted disk by property~(3) 
in Proposition~\ref{thm1}.

Finally we consider the case where $L(\m)$ is a Hopf multilink.
Let $m_1S_1$ and $m_2S_2$ denote the multilink components of $L(\m)$, i.e.,
$L(\m)=(\Sigma(1,1), m_1S_1\cup m_2S_2)$.
If $m_1+m_2\ne 0$ then $L(\m)$ satisfies (PTP) up to the reversal of the orientation of $L(\m)$.
So, the above proof works in this case.
Suppose that $m_1+m_2=0$. Since the orientation of $L(\m)$ is not canonical, it is enough to
check that the compatible contact structure is overtwisted.
This follows immediately since
the fiber surface of $L(\m)$ is a disjoint union of the fiber surfaces of a negative Hopf link
and the compatible contact structure is same as that of the negative Hopf link.
\qed
\vspace{3mm}

\section{Case $a_1a_2\cdots a_k<0$}

\subsection{Explicit construction of the contact structure}

Throughout this section, we assume that $A=a_1\cdots a_k<0$.
We start from the following lemma.

\begin{lemma}\label{lemma20}
If $A<0$ then every fibered Seifert multilink has at least one negative component.
\end{lemma}

\begin{figure}[htbp]
   \centerline{\input{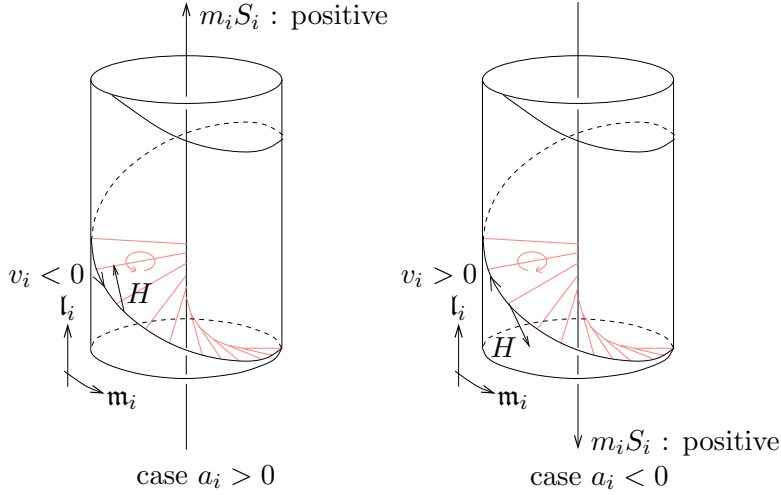}}
   \caption{The framing of $F$ along $m_iS_i$.\label{fig6}}
\end{figure}

\begin{proof}
The proof is analogous to that of Lemma~\ref{lemma19}.
In the present case, the framing of the fiber surface $F$ along $m_iS_i$ becomes
as shown in Figure~\ref{fig6}, from which we have $a_iv_i<0$.
Hence the right hand side of equation~\eqref{eq100} is
strictly negative since $a_iu_i-\sigma_iv_i>0$ and $a_iv_i<0$.
This is a contradiction.
\end{proof}

The main assertion in this section is the following.

\begin{prop}\label{thm2}
Let $L(\m)=(\Sigma, m_1S_1\cup\cdots\cup m_nS_n)$
be a fibered Seifert multilink $L(\m)$ in a homology $3$-sphere $\Sigma=\Sigma(a_1,\ldots,a_k)$ with $A<0$.
Assume {\rm (PTP)}.
Fix an index $i_0$ of some negative component of $L(\m)$.
Then there exists a positive contact form $\alpha$ on $\Sigma$
with the following properties:
\begin{itemize}
\item[(1)] $L(\m)$ is compatible with the contact structure $\xi=\ker\alpha$.
\item[(2)] The Reeb vector field $R_\alpha$ of $\alpha$ 
is tangent to the fibers of the Seifert fibration on $\mathcal S\times S^1$.
\item[(3)] The neighborhood $(D^2\times S^1)_i$ of each negative component $m_iS_i$, 
except $m_{i_0}S_{i_0}$, 
contains a Lutz tube. In particular, it contains an overtwisted disk.
\item[(4)] On the other $(D^2\times S^1)_i$'s, 
except $i=i_0$, $\ker\alpha$ is transverse to the fibers of the Seifert fibration.
\end{itemize}
In particular, if $L(\m)$ has at least two negative components then
the contact structure $\ker\alpha$ is overtwisted.
\end{prop}

Before proving this proposition, we prepare a lemma similar to Lemma~\ref{lemma1}.

\begin{lemma}\label{lemma2}
Suppose $A<0$ and fix an index $i_0$.
Let $U_i$ be a collar neighborhood of $C_i$ in $\mathcal S$ 
with coordinates $(r_i,\theta_i)\in [1,2)\times S^1$ 
satisfying $\{(r_i,\theta_i) \,;\, r_i=1\}=C_i$.
Then there exists a $1$-form $\beta$ on $\mathcal S$ which satisfies the following properties:
\begin{itemize}
\item[(1)] $d\beta>0$ on $\mathcal S$.
\item[(2)] If $b_i/a_i\leq 0$ and $i\ne i_0$ then 
$\beta=R_ir_id\theta_i$ with $-b_i/a_i<R_i$ near $C_i$ on $U_i$.
\item[(3)] If $b_i/a_i>0$ and $i\ne i_0$ then
$\beta=(R_i/r_i)d\theta_i$ with $-b_i/a_i<R_i<0$ near $C_i$ on $U_i$.
\item[(4)] If $b_{i_0}/a_{i_0}-1/A < 0$ then 
$\beta=R_{i_0}r_{i_0}d\theta_{i_0}$ with $0<R_{i_0}<-b_{i_0}/a_{i_0}+1/A$ 
 near $C_{i_0}$ on $U_{i_0}$.
\item[(5)] If $b_{i_0}/a_{i_0}-1/A\geq 0$ then
$\beta=(R_{i_0}/r_{i_0})d\theta_{i_0}$ 
with $R_{i_0}<-b_{i_0}/a_{i_0}+1/A$  near $C_{i_0}$ on $U_{i_0}$.
\end{itemize}
\end{lemma}

\begin{proof}
Since $\sum_{i\ne i_0}(-b_i/a_i)+(-b_{i_0}/a_{i_0}+1/A)=0$,
we can choose $R_1,\ldots,R_k$ such that they satisfy the above inequalities and 
the inequality $\sum_{i=1}^k R_i<0$.
The $1$-form $\beta$ required can be constructed from these $R_i$'s
in the same way as in the proof of Lemma~\ref{lemma1}.
\end{proof}

\noindent
{\it Proof of Proposition}~\ref{thm2}.
We make a contact form $\alpha_0$ on $\mathcal S\times S^1$ from the $1$-form $\beta$
in Lemma~\ref{lemma2}
and extend it to $(D^2\times S^1)_i$ as in the proof of Proposition~\ref{thm1}.
Properties~(2), (3), (4) in the assertion follow from this construction. 
Let $\alpha$ denote the obtained contact form on $M$.

\begin{figure}[htbp]
   \centerline{\input{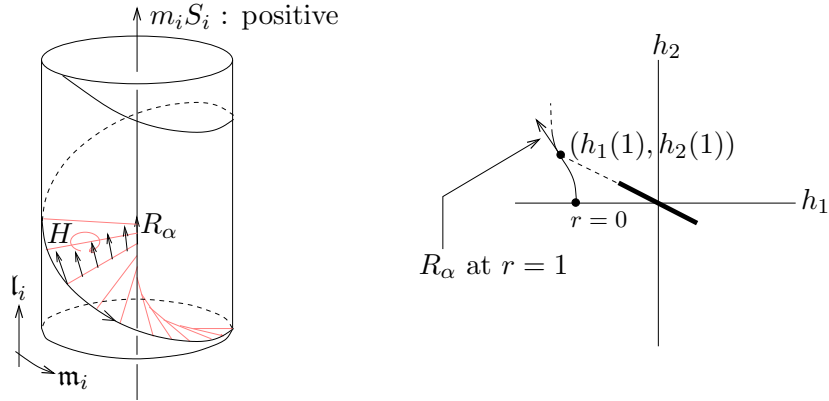}}
   \caption{The mutual positions of $F$, $H$ and $R_\alpha$ in the case where
$m_iS_i$ is a positive component.\label{fig17}}
\end{figure}

Suppose that $i\ne i_0$, $m_iS_i$ is a positive component and $a_i>0$.
From equation~\eqref{eqaaa}, 
we have $h_1(1)<0$, $h_2(1)>0$, $h'_1(1)<0$ and $h'_2(1)>0$.
Hence the mutual positions of the fiber surface $F$, the oriented fibers $H$ of the Seifert
fibration and the Reeb vector field $R_\alpha$ on $(D^2\times S^1)_i$ 
become as shown on the left in Figure~\ref{fig17}.
The contact structure $\alpha$ in this case is determined by the curve described on
the right. 
From these figures, we can easily check that these satisfy property~(1) in the assertion.

\begin{figure}[htbp]
   \centerline{\input{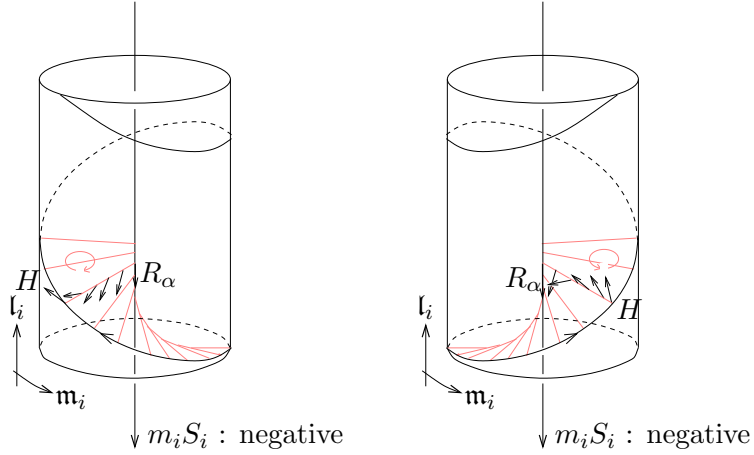}}
   \caption{The mutual positions of $F$, $H$ and $R_\alpha$ in the case where
$m_iS_i$ is a negative component.\label{fig18}}
\end{figure}

If $m_iS_i$ is negative and $a_i>0$ then we have the same inequalities.
Hence their mutual positions become as shown in Figure~\ref{fig18} and the property~(1) holds.
If $i=i_0$ then $h_2(1)>0$ may not hold, but this does not make any problem since
$m_{i_0}S_{i_0}$ is a negative component.
Thus the property~(1) holds.

The proof is analogous in case $a_i<0$.
\qed

\subsection{Some criterion to detect overtwisted disks}

In this subsection, we show two lemmas which give sufficient conditions for 
the contact structure in Proposition~\ref{thm2} to be overtwisted.

\begin{lemma}\label{lemma3}
Suppose $A<0$ and let $m_{i_0}S_{i_0}$ be a negative component of $L(\m)$.
Suppose further that there exists $a_{i_1}$ among $a_1,\ldots,a_k$ which satisfies the inequality
\[
   \frac{1}{|a_{i_1}|}\left(\frac{1}{|a_{i_0}|}-\frac{1}{|a_{i_1}|}\right)>-\frac{1}{A}.
\]
Then the contact structure in Proposition~\ref{thm2} is overtwisted.
\end{lemma}

\begin{proof}
From the inequality in the assumption,
we have $|a_{i_1}|>|a_{i_0}|$. In particular, $i_0\ne i_1$.
We can assume that $m_{i_1}S_{i_1}$ is a positive component, since otherwise
the contact structure is overtwisted by Proposition~\ref{thm2}.
We will find $R_1,\ldots,R_k$ in Lemma~\ref{lemma2} which satisfy
\[
   |a_{i_0}|\left(R_{i_0}+\frac{b_{i_0}}{a_{i_0}}\right)
   =-|a_{i_1}|\left(R_{i_1}+\frac{b_{i_1}}{a_{i_1}}\right)<0.
\]
Set $X=R_{i_0}+b_{i_0}/a_{i_0}$ and $Y=R_{i_1}+b_{i_1}/a_{i_1}$.
They should satisfy the conditions in Lemma~\ref{lemma2}, that is,
$X-1/A<0$ and $Y>0$.

For a sufficiently small $\ve>0$, we set $R_i$'s for $i\ne i_0, i_1$ 
such that they satisfy the conditions in Lemma~\ref{lemma2} and the equality
\[
   \sum_{i\ne i_0, i_1}\left(R_i+\frac{b_i}{a_i}\right)=\ve.
\]
In the case $k=2$, we set $\ve=0$.
We need the inequality $\sum_{i=1}^kR_i<0$ and hence $X$ and $Y$ should satisfy
\[
   0>\sum_{i\ne i_0,i_1}R_i + R_{i_0} + R_{i_1}
=\ve-\sum_{i\ne i_0,i_1}\frac{b_i}{a_i}+ R_{i_0} + R_{i_1}
=\ve-\frac{1}{A}+X+Y.
\]

Now we assume that the following inequality holds:
\begin{equation}\label{eq1000}
   |b_{i_0}+a_{i_0}R_{i_0}|=-|a_{i_0}|X<\frac{1}{|a_{i_0}|}.
\end{equation}
The difference of the slopes of a meridional disk and a Legendrian curve
on $\bd (D^2\times S^1)_{i_0}$ is given as
\[
   \frac{(a_{i_0}Q_{i_0}+b_{i_0}H)}{a_{i_0}}-(Q_{i_0}-R_{i_0}H)
=\left(\frac{b_{i_0}}{a_{i_0}}+R_{i_0}\right)H.
\]
Since $b_{i_0}/a_{i_0}+R_{i_0}=X<1/A<0$, the slope of the Legendrian curve
is a bit higher than that of the meridional disk, see Figure~\ref{fig7}.
Let $\gamma$ be the boundary of the meridional disk.
Since the distance of two neighboring intersection points of $H$ and $\gamma$ is $1/|a_{i_0}|$,
inequality~\eqref{eq1000} ensures that
we can isotope $\gamma$ on $\bd (D^2\times S^1)_{i_0}$ such that
it is Legendrian except for a short vertical interval of length $|b_{i_0}+a_{i_0}R_{i_0}|$.
We denote by $\Delta_{i_0}$ the meridional disk bounded by this isotoped $\gamma$.
\begin{figure}[htbp]
   \centerline{\input{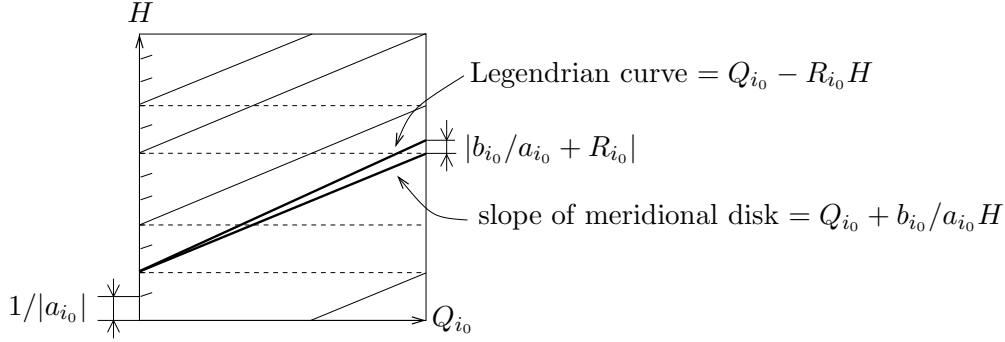}}
   \caption{The slopes of a meridional disk and a Legendrian curve
on the boundary of $(D^2\times S^1)_{i_0}$.\label{fig7}}
\end{figure}

We also obtain a similar disk $\Delta_{i_1}$ in $(D^2\times S^1)_{i_1}$,
assuming the inequality
\[
   |b_{i_1}+a_{i_1}R_{i_1}|=|a_{i_1}|Y<\frac{1}{|a_{i_1}|}.
\]
In this case, the slope of the Legendrian curve
is a bit lower than that of the meridional disk since 
$b_{i_1}/a_{i_1}+R_{i_1}=Y>0$, cf.~Figure~\ref{fig9}.

In summary, we have assumed for a point $(X,Y)$ to satisfy the following conditions:
\begin{equation}\label{eq1001}
\left\{
\begin{split}
& |a_{i_0}|X+|a_{i_1}|Y=0,  \\
& X+Y<-\ve+\frac{1}{A}, \\
& -\frac{1}{a_{i_0}^2}<X<\frac{1}{A}, \\
& 0<Y<\frac{1}{a_{i_1}^2}.
\end{split}
\right.
\end{equation}
Note that we always have the inequality $-1/a_{i_0}^2<1/A$,
because $1/|a_{i_1}|(1/|a_{i_0}|-1/|a_{i_1}|)>-1/A$ implies
 $|a_{i_0}|<|a_{i_1}|$ and hence
\[
  -\frac{1}{a_{i_0}^2}<-\frac{1}{|a_{i_0}||a_{i_1}|}\leq \frac{1}{A}.
\]

Now we describe the region on the $XY$-plane where $(X,Y)$ satisfies 
the inequalities in the above conditions, which is shown in Figure~\ref{fig8}.
\begin{figure}[htbp]
   \centerline{\input{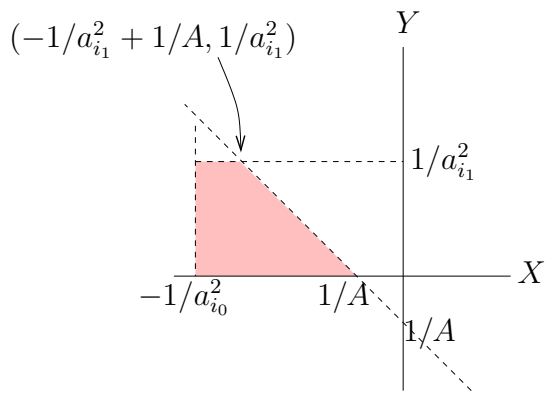}}
   \caption{The region where $(X,Y)$ satisfies the required inequalities.\label{fig8}}
\end{figure}
Note that we used the inequality 
\[
   \frac{1}{a_{i_0}^2}-\frac{1}{a_{i_1}^2}
   >\frac{1}{|a_{i_1}|}\left(\frac{1}{|a_{i_0}|}-\frac{1}{|a_{i_1}|}\right)>-\frac{1}{A}
\]
when we described this region.
The equality and inequalities in~\eqref{eq1001} have a solution if and only if
the line $|a_{i_0}|X+|a_{i_1}|Y=0$ intersects this region, i.e.,
the following inequality holds:
\[
   |a_{i_0}|\left(-\frac{1}{a_{i_1}^2}+\frac{1}{A}\right)+|a_{i_1}|\left(\frac{1}{a_{i_1}^2}\right)>0,
\]
and this follows from the assumption.
Thus the embedded disks $\Delta_{i_0}$ and $\Delta_{i_1}$ exist.

Finally we connect these disks by a band $B$ whose two sides are Legendrian
as shown in Figure~\ref{fig9}. We here explain this more precisely.
We first remark that 
the lengths of the two short vertical intervals on the boundaries of $\Delta_{i_0}$ and $\Delta_{i_1}$
are the same since
\[
   |b_{i_0}+a_{i_0}R_{i_0}|=-|a_{i_0}|X=|a_{i_1}Y=|b_{i_1}+a_{i_1}R_{i_1}|.
\]
Let $p_0$, $q_0$ be the endpoints of the vertical interval of the boundary of $\Delta_{i_0}$ 
and let $p_1$ and $q_1$ be those of $\Delta_{i_1}$.
Choose a vertical annulus $W=H\times [0,1]$ between $(D^2\times S^1)_{i_0}$ and $(D^2\times S^1)_{i_1}$
as shown in Figure~\ref{fig9}
and let $\mathcal F_W$ denote the foliation on $W$ determined by $\xi$.
Note that $\mathcal F_W$ is non-singular and every leaf of $\mathcal F_W$ connects
the connected components of $\bd W$ because $\xi$ is transverse to $H$.
By shifting $\Delta_{i_0}$ if necessary, we can assume that
$p_0$ and $p_1$ are the endpoints of the same leaf of $\mathcal F$.
Since the lengths of the short vertical intervals are the same,
by shifting both of $\Delta_{i_0}$ and $\Delta_{i_1}$ simultaneously,
we can find positions of $\Delta_{i_0}$ and $\Delta_{i_1}$ such that
$p_0$ and $p_1$ are the endpoints of a leaf of $\mathcal F$ and 
$q_0$ and $q_1$ are also the endpoints of another leaf of $\mathcal F$.
Now we choose the band $B$ to be a curved rectangle
such that its boundary consists of these leaves and the short vertical intervals
and it is tangent to the contact structure $\xi$ along the leaves of $\mathcal F_W$ on the boundary
as shown in Figure~\ref{fig9}.
The union $\Delta_{i_0}\cup B\cup \Delta_{i_1}$ is a disk embedded in $\Sigma$ with
polygonal Legendrian boundary.
We then isotope it in a neighborhood of the corners of the polygonal Legendrian boundary 
such that it becomes a smooth embedded disk with smooth Legendrian boundary.
From the construction, the contact structure $\xi$ is tangent to this disk
along its boundary. Hence it is an overtwisted disk.
\end{proof}

\begin{figure}[htbp]
   \centerline{\input{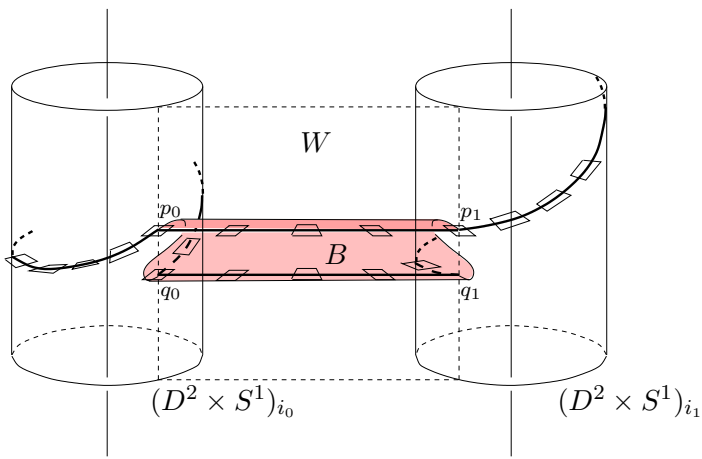}}
   \caption{The band $B$.\label{fig9}}
\end{figure}

\begin{lemma}\label{lemma4}
Suppose $A<0$ and let $m_{i_0}S_{i_0}$ be a negative component of $L(\m)$.
Suppose further that there exist $a_{i_1}$ and $a_{i_2}$ satisfying $|a_{i_0}|<|a_{i_2}|<|a_{i_1}|$.
Then the contact structure in Proposition~{\rm \ref{thm2}} is overtwisted.
\end{lemma}

\begin{proof}
We have the inequality
\[
   -\frac{|a_{i_1}|}{A}\leq\frac{1}{|a_{i_0}a_{i_2}|}
=\left(\frac{1}{|a_{i_0}|}-\frac{1}{|a_{i_2}|}\right)\frac{1}{|a_{i_2}|-|a_{i_0}|}
\leq \frac{1}{|a_{i_0}|}-\frac{1}{|a_{i_2}|}<\frac{1}{|a_{i_0}|}-\frac{1}{|a_{i_1}|}
\]
and hence the assertion follows from Lemma~\ref{lemma3}.
\end{proof}

\begin{ex}
Suppose that $\gcd(|p|,|q|)=1$ and $pq<0$.
\begin{itemize}
\item[(1)] 
$(\Sigma,L)=(\Sigma(1,p,q),-S_1)$ is a $(p,q)$-torus knot in $S^3$. Here
the component $-S_1$ must be negative because of Lemma~\ref{lemma20}.
If $|p|, |q|\geq 2$ then there exists an overtwisted disk by Lemma~\ref{lemma4}.
If either $|p|=1$ or $|q|=1$ then $L$ is a trivial knot in $S^3$ and
its compatible contact structure is tight.
Actually, this does not satisfy the condition in Lemma~\ref{lemma3}.
\item[(2)] 
$(\Sigma,L)=(\Sigma(p,q),S_1\cup -S_2)$ is a positive Hopf link in $S^3$.
It is well-known that its compatible contact structure is tight, and
this actually does not satisfy the condition in Lemma~\ref{lemma3}.
\end{itemize}
\end{ex}

\section{Fibered Seifert links in $S^3$}

In this section, we study Seifert links in $S^3$.
The classification of Seifert links in $S^3$ was done by Burde and Murasugi~\cite{bm},
in which they proved that
a link is a Seifert link in $S^3$ if and only if it is a union of a finite number of fibers
of the Seifert fibration in $\Sigma(p,q)$ with $pq\ne 0$ or $(p,q)=(0,1)$ (cf.~\cite[p.62]{en}).
The classification of contact structures on $S^3$ had been done by Eliashberg~\cite{eliashberg:89,eliashberg}.
In particular, it is known that $S^3$ admits a unique tight contact structure
up to contactomorphism, so-called the {\it standard contact structure}.
\vspace{3mm}

\noindent
{\it Proof of Theorem}~\ref{cor03}.
The assertion in case $pq>0$ follows from Theorem~\ref{thm01}.
Suppose $pq<0$.
We first prove the assertion in the case where all components of $L$ are negative.
In this case, (PTP) is satisfied by Lemma~\ref{lemma20}.
If $L$ has more than one link components then
the contact structure is overtwisted by the last assertion in Proposition~\ref{thm2}.
Suppose that $L$ consists of only one component, then $L$ is either a trivial knot or 
a $(p,q)$-torus knot with $pq<0$.
It is well-known that the contact structure of a trivial knot is tight, and that
the contact structure of a $(p,q)$-torus knot with $pq<0$ is 
overtwisted if and only if it is not a trivial knot. Thus the assertion follows in this case.

Next we consider the case where $L$ has at least one positive component.
Note that $L$ also has one negative component by Lemma~\ref{lemma20}.
We can assume that the number of negative components of $L$ is one, otherwise
the contact structure is overtwisted by the last assertion in Proposition~\ref{thm2}.

We decompose the argument into three cases:
\begin{itemize}
\item[(1)] The two exceptional fibers of $\Sigma(p,q)$ are both components of $L$. That is,
\[
   L=(\Sigma(\underbrace{1, \ldots, 1}_{n-2},p,q),
m_1S_1\cup\cdots\cup m_{n-2}S_{n-2}\cup m_{n-1}S_{n-1}\cup m_nS_n).
\]
\item[(2)] One of the two exceptional fibers of $\Sigma(p,q)$ is a component of $L$. That is,
\[
   L=(\Sigma(\underbrace{1, \ldots, 1}_{n-1},p,q), m_1S_1\cup\cdots\cup m_{n-1}S_{n-1}\cup m_nS_n).
\]
\item[(3)] Neither of the two exceptional fibers of $\Sigma(p,q)$ is a component of $L$. That is,
\[
   L=(\Sigma(\underbrace{1, \ldots, 1}_n,p,q),m_1S_1\cup\cdots\cup m_nS_n).
\]
\end{itemize}
Here $m_i\in\{-1,+1\}$ since $L$ is a fibered link.

\begin{figure}[htbp]
   \centerline{\input{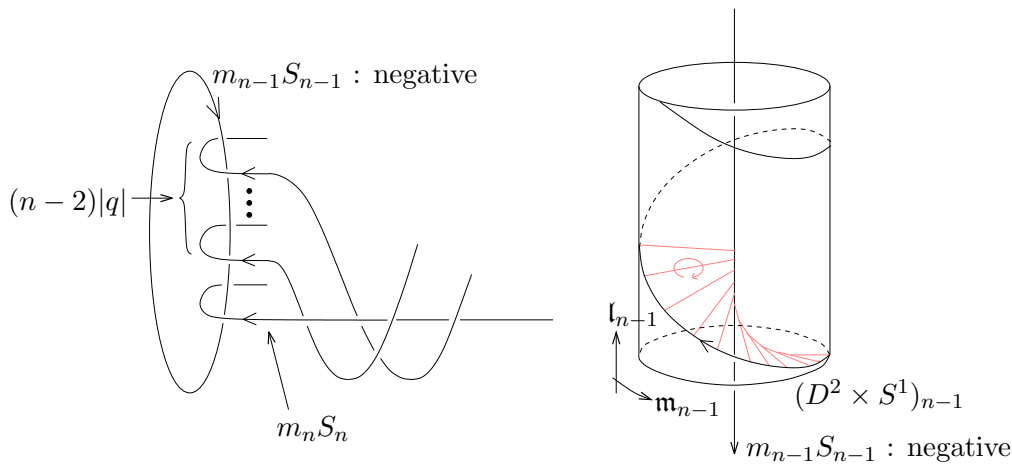}}
   \caption{The framing of the Seifert surface in case (1) with 
negative component $m_{n-1}S_{n-1}$ and $p>0$.\label{fig10}}
\end{figure}

We first consider case~(1).
If $n=2$ then $L$ is a positive Hopf link in $S^3$.
Suppose $n\geq 3$ and that either $S_{n-1}$ or $S_n$, say $S_{n-1}$, is a negative component.
The linking number of $m_{n-1}S_{n-1}$ and all the other components of $L$ is $(n-2)|q|+1$.
Note that $n-2$ is the number of the link components of $L$ along non-exceptional fibers.
For a fiber surface $F$ of $L$,
the oriented boundary $\bd(F\cap (D^2\times S^1)_{n-1})\setminus m_{n-1}S_{n-1}$
on $\bd (D^2\times S^1)_{n-1}$
is given as $\gamma=\pm (-((n-2)|q|+1)\fM_{n-1}+\fL_{n-1})$, where the sign $\pm$ is $+$ if
$p>0$ and $-$ otherwise, see Figure~\ref{fig10}.
Here the surface on the right is described by applying the Seifert's algorithm
to the diagram on the left.

Since $H=q\fM_{n-1}+p\fL_{n-1}$, (PTP) implies the inequality 
$I(\gamma,H)=\mp(((n-2)|q|+1)p+q)>0$, where $I(\gamma,H)$ is the algebraic intersection number
of $\gamma$ and $H$ on $\bd (D^2\times S^1)_{n-1}$.
However, 
\[
\begin{split}
   I(\gamma,H)&=\mp(((n-2)|q|+1)p+q)=(n-2)pq\mp(p+q)\\
&=(p\mp 1)(q\mp 1)+(n-3)pq-1<0
\end{split}
\]
since $(p\mp 1)(q\mp 1)\leq 0$ and $(n-3)pq\leq 0$ for $n\geq 3$.
This is a contradiction.

\begin{figure}[htbp]
   \centerline{\input{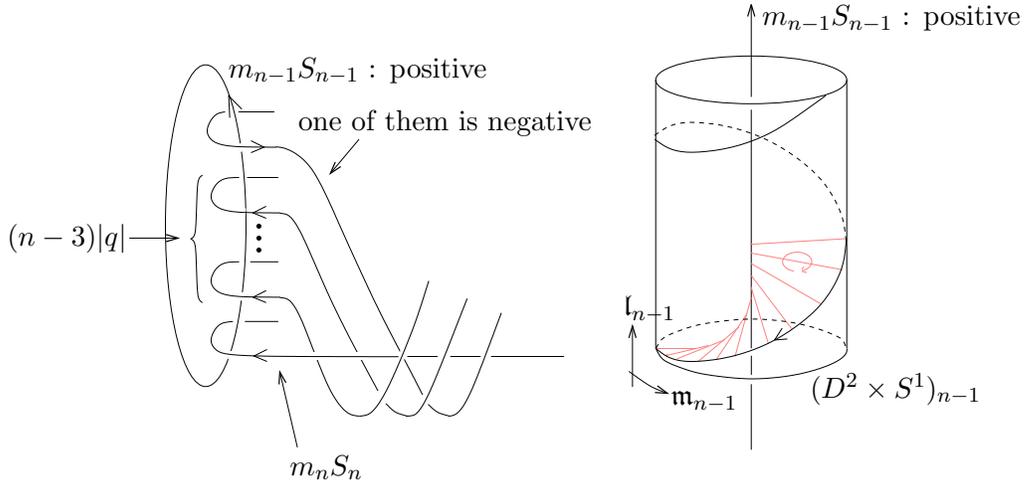}}
   \caption{The framing of the Seifert surface in case~(1) with a non-exceptional
fiber being the negative component.\label{fig11}}
\end{figure}

Suppose $n\geq 3$ and a regular fiber is a negative component of $L$.
The linking number of $m_{n-1}S_{n-1}$ and
all the other components of $L$ is $-(n-4)|q|-1$ and 
the oriented boundary $\bd(F\cap (D^2\times S^1)_{n-1})\setminus m_{n-1}S_{n-1}$
on $\bd (D^2\times S^1)_{n-1}$
becomes $\gamma=\pm ((-(n-4)|q|-1)\fM_{n-1}-\fL_{n-1})$, see Figure~\ref{fig11}.
Thus, $I(\gamma,H)=\mp(((n-4)|q|+1)p-q)=(n-4)pq\mp p\pm q$.
If $|p|,|q|\geq 2$ then the contact structure of $L$ is overtwisted by Lemma~\ref{lemma4}.
If either $|p|$ or $|q|$ equals $1$ then
\[
(n-4)pq\mp p\pm q=(n-3)pq-(p\mp 1)(q\pm 1)-1<0
\]
since $(p\mp 1)(q\pm 1)=0$. Hence (PTP) does not hold.

Next we consider case~(2).
If $n=1$ then $L$ is a trivial knot in $S^3$. 
Suppose $n\geq 2$ and that $S_n$ is a negative component.
Since
\[
   I(\gamma,H)=\mp((n-1)|q|p+q)=(n-1)pq\mp q=(n-1)pq+|q|\leq 0,
\]
(PTP) does not hold (cf.~Figure~\ref{fig10} with deleting the component $m_nS_n$ and
replacing the number $(n-2)|q|$ by $(n-1)|q|$ and the indices $n-1$ by $n$).
We remark that the equality holds when $n=2$ and $|p|=1$, and if $|q|=1$ in addition then
$L$ becomes a positive Hopf link. Nevertheless,
we can ignore this case because the fibration of a positive
Hopf link is not given by this Seifert fibration.

Suppose $n\geq 2$ and a regular fiber is a negative component of $L$,
then
\[
I(\gamma,H)=\mp((n-3)|q|p-q)=(n-3)pq\pm q=(n-3)pq-|q|
\]
(cf.~Figure~\ref{fig11} with deleting the component $m_nS_n$
and replacing the number $(n-3)|q|$ by $(n-2)|q|$ and the indices $n-1$ by $n$).
This is positive if and only if $n=2$ and $|p|\geq 2$, in which case
if $|q|\geq 2$ then the contact structure of $L$ is overtwisted by Lemma~\ref{lemma4},
and if $|q|=1$ then $L$ is a positive Hopf link and its contact structure is tight.

\begin{figure}[htbp]
   \centerline{\input{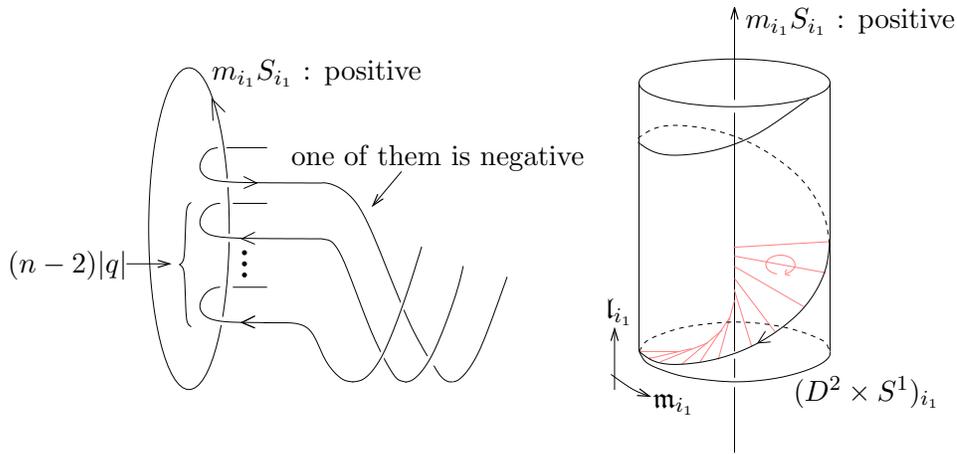}}
   \caption{The framing of the Seifert surface in case (3).\label{fig12}}
\end{figure}

Finally we consider case~(3).
If $n=1$ then it is a $(p,q)$-torus knot and we know that
its contact structure is tight if and only if it is a trivial knot.
If $n=2$ then $L$ is a positive Hopf link, otherwise $L$ is not fibered.
If $n\geq 3$ and $|p|, |q|\geq 2$ then its contact structure is overtwisted
by Lemma~\ref{lemma4}. So, we can suppose that $n\geq 3$ and either $|p|$ or $|q|$ equals $1$.
Choose a positive component $m_{i_1}S_{i_1}$ of $L$, then 
the oriented boundary $\bd(F\cap (D^2\times S^1)_{i_1})\setminus m_{i_1}S_{i_1}$
on $\bd (D^2\times S^1)_{i_1}$
is given as $\gamma=-(n-3)|q|\fM_{i_1}-\fL_{i_1}$, see Figure~\ref{fig12}.
Since $I(\gamma,H)=-(n-3)|q|+pq<0$, (PTP) does not hold.

If $pq=0$ then $L$ is as shown in Figure~\ref{fig111},
which is a connected sum of a finite number of Hopf links.
The plumbing argument in~\cite{torisu} ensures that
the contact structure of such a link is tight
if and only if every summand is a positive Hopf link. This completes the proof.
\qed

\section{Seifert links in $S^3$ and their strongly quasipositivity}

A Seifert surface in $S^3$ is called {\it quasipositive} if 
it is obtained from a finite number of parallel copies of a disk 
by attaching positive bands.
A link is called {\it strongly quasipositive} if it is realized as
the boundary of some quasipositive surface.
In other words, a strongly quasipositive link is the closure of a braid
given by the product of words of the form
\[
   \sigma_{i,j}=(\sigma_i\cdots\sigma_{j-2})\sigma_{j-1}(\sigma_i\cdots\sigma_{j-2})^{-1}
\]
where $\sigma_i$ is a positive generator of braid.
See \cite{rudolph1,rudolph2,rudolph3,rudolph4,rudolph5,rudolph6}
for further studies of quasipositive surfaces.

It is known by Hedden~\cite{hedden}, and Baader and the author~\cite{bi} 
in a different way, that the compatible contact structure of a fibered link in $S^3$ is tight
if and only if its fiber surface is quasipositive.
So, Theorem~\ref{cor03} can be generalized into the non-fibered case
as stated in Corollary~\ref{cor04}.
\vspace{3mm}

\noindent
{\it Proof of Corollary}~\ref{cor04}.
The assertion had been proved in Theorem~\ref{cor03} if $L$ is fibered.
So, hereafter we assume that $L$ is non-fibered.
If $(a_1,a_2)=(0,1)$ then $L$ must be a trivial link with several components,
which is excluded by the assumption.
Suppose that $a_1a_2\ne 0$.
By using the criterion in~\cite[Theorem~11.2]{en}, we can easily check that
$L$ is not fibered if and only if it is a positive or negative torus link,
other than a Hopf link,
which consists of even number of link components, say $2k$, half of which
have reversed orientation.
Such an $L$ is realized as the boundary of a Seifert surface $F$
consisting of $k$ annuli.

Suppose $a_1a_2>0$ and let $F'$ be one of the annuli of $F$.
The core curve of $F'$ constitutes a positive torus knot,
say a $(p,q)$ torus knot with $p, q>0$.
It is known in~\cite[Lemma~6.1]{bi} that if $F'$ is quasipositive then
$-1$ times the linking number $lk(F')$ of the two boundary components of $F'$
is at most the maximal Thurston-Bennequin number $TB(K)$
of the core curve $K$ of the annulus, i.e. $-lk(F')\leq TB(K)$.
It is known in~\cite{tanaka} that
\[
   TB(K)=(p-1)q-p=pq-p-q,
\]
where we regarded $p$ as the number of Seifert circles, which equals the braid index.
However, we can easily check $lk(F')=-pq$, which does not satisfy the inequality $-lk(F')\leq TB(K)$.
Thus $F'$ is not quasipositive.
Now assume that $L$ is strongly quasipositive. Then, by definition,
there exists a quasipositive surface bounded by $L$.
However this surface contains the above non-quasipositive annulus as an essential subsurface,
which contradicts the Characterization Theorem of quasipositive surfaces in~\cite{rudolph1}.
Thus $L$ is not strongly quasipositive.

\begin{figure}[htbp]
   \centerline{\input{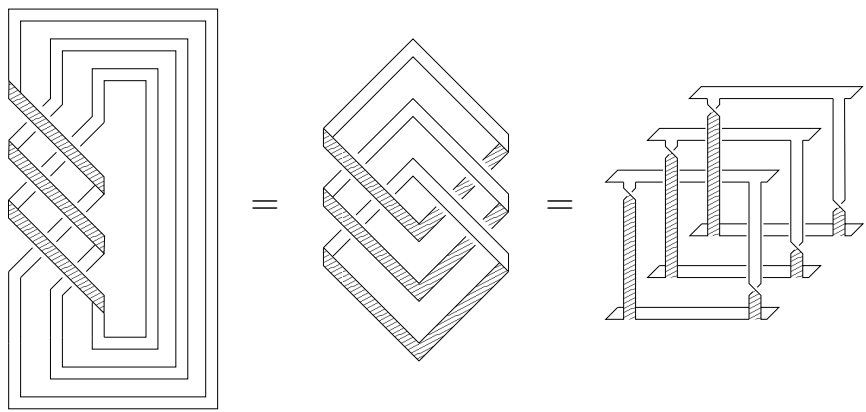}}
   \caption{The surface $F$ in the case $(p,q)=(3,-3)$.\label{fig112}}
\end{figure}

\begin{figure}[htbp]
   \centerline{\input{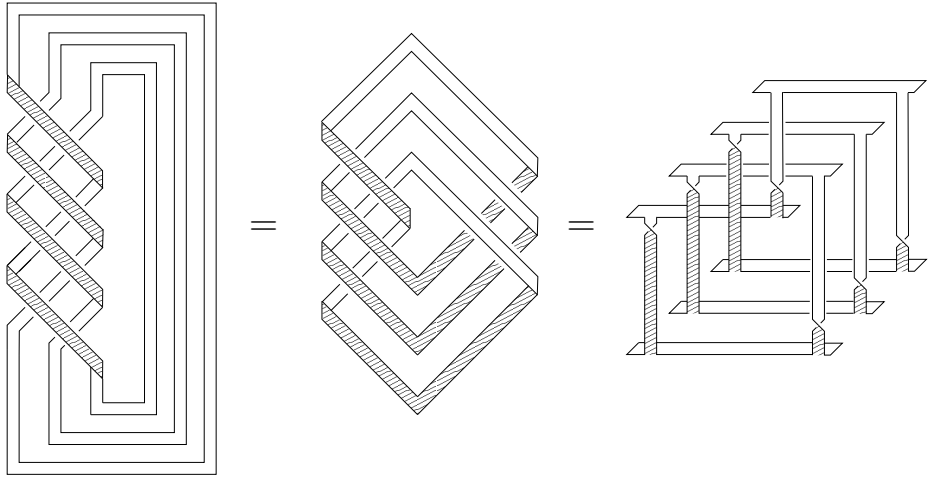}}
   \caption{The surface $F$ in the case $(p,q)=(3,-4)$.\label{fig113}}
\end{figure}

If $a_1a_2<0$ then the link $L$ is in case~(3) in the assertion.
Suppose that the core curves of annuli of $F$ constitutes a $(kp,kq)$ torus link 
with $p>0$ and $q<0$.
Using ambient isotopy move in $S^3$, we can assume that $p\leq |q|$.
In the case where $p=|q|$, 
we set the surface $F$ in the position as shown in Figure~\ref{fig112},
which shows that the surface is quasipositive.
If $p<|q|$, we need to add more crossings, though we can check
that the surface is still quasipositive  as shown in Figure~\ref{fig113}.
This completes the proof.
\qed

\vspace{3mm}

We close this section with a conjecture arising from the fact in Corollary~\ref{cor04}.

\begin{conj}
Any non-splittable unoriented link in $S^3$ has at most two strongly quasipositive orientations.
\end{conj}

Here a strongly quasipositive orientation means an orientation assigned to
the unoriented link such that the obtained oriented link
becomes strongly quasipositive.
As in Corollary~\ref{cor04}, this conjecture is true for all Seifert links in $S^3$.
We will prove the same assertion for fibered, positively-twisted graph links
in $S^3$ in the subsequent paper~\cite{ishikawa3}.

\section{Cablings}

\subsection{Definition of positive and negative cablings}

In this section, we study a fibered multilink in a $3$-manifold with cabling structures.
The notion of multilink is convenient to describe relation between compatible contact structures
before and after the cabling. For this aim, we will give a definition of cabling in an unusual way.
Our definition coincides with the usual definition of cabling in the case
where the cabling is performed along a fibered knot in a $3$-manifold.
This will be discussed in Corollary~\ref{cor7}.

Let $M$ be an oriented, closed, smooth $3$-manifold and $L(\m)$ a fibered multilink in $M$.
Suppose that there exists a solid torus $N$ in $M$ such that each $L(\m)\cap N$ is a torus
multilink in $N$ with consistent orientation, 
i.e., a multilink in $N$ lying on a torus parallel to the boundary $\bd N$ 
all of whose link components have consistent orientations.
We replace the torus multilink component of $L(\m)$ in $N$ by its core curve $S$,
extend the fiber surfaces of $L(\m)$ by the retraction of $N$ to $S$,
and define the multiplicity of $S$ from these fiber surfaces canonically.
We denote the obtained multilink in $M$ by $L'(\m')$. Note that $L'(\m')$ is always fibered.
The operation producing $L(\m)$ from $L'(\m')$ by attaching $L(\m)\cap N$ along $S$
is called a {\it cabling}.

Next we define the notion of positive and negative cablings.
We set $L(\m)\cap N$ and $L'(\m')$ in $M$ simultaneously such that the core curve of $N$ coincides with 
the link component of $L'(\m')$ in $N$, and check the intersection of 
$L(\m)\cap N$ with the fiber surface of $L'(\m')$.
Note that this intersection is always transverse, see Lemma~\ref{lemma3001}  below.

\begin{dfn}\label{dfn04}
A cabling is called {\it positive} if $L(\m)\cap N$ 
intersects the fiber surface of $L'(\m')$  positively transversely.
If the intersection is negative then the cabling is called {\it negative}.
\end{dfn}

To discuss the framing of the cabling, we fix a basis of $\bd N$ as follows:
Let $\fM$ be an oriented meridian on $\bd N$ positively transverse to the fiber surface
$F$ of $L(\m)$ and 
$\fL$ be an oriented simple closed curve on $\bd N$ such that $I(\fM,\fL)=1$,
where $I(\fM,\fL)$ is the algebraic intersection number of $\fM$ and $\fL$ on $\bd N$.
Each connected component of the oriented boundary of $F\setminus \text{int}N$ 
on $\bd (M\setminus\text{int\,}N)$ is given as
$\gamma=u\fM+v\fL$, where $(u,v)\in\Z\times\N$ are assumed to be coprime.

\begin{figure}[htbp]
   \centerline{\input{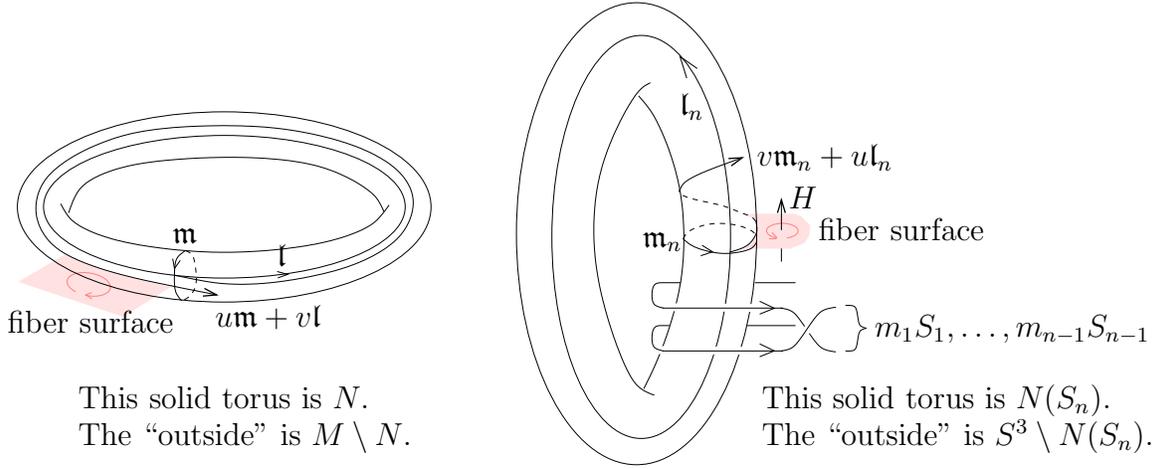}}
   \caption{The left figure shows the fiber surface $F$ in $M\setminus\text{int\,} N$ 
and the right one shows $L(\m)\cap N$ in $N\subset S^3$.\label{fig13}}
\end{figure}

Now we embed $N$ into $S^3$ along a trivial knot such that
$(\fM,\fL)$ becomes the preferred meridian-longitude pair of this trivial knot.
We then add the core curve $S_n$ of $S^3\setminus\text{int\,}N$ as an additional link component
to $L(\m)\cap N$ embedded in $S^3$, extend the fiber surfaces of $L(\m)$ by the retraction of 
$S^3\setminus\text{int\,}N$ to $S_n$,
and define the multiplicity $m_n$ of $S_n$ from these fiber surfaces canonically.
The obtained multilink can be represented as
\[
   L_{p,q}(\m_{p,q})=(\Sigma(\underbrace{1,\ldots,1}_{n-1},\ve q,\ve\,p),
 m_1S_1\cup\cdots\cup m_{n-1}S_{n-1}\cup \ve_n m_nS_n),
\]
where $p>0$,
\[
\ve=
   \begin{cases}
   1 & \text{if the cabling is positive}  \\
   -1 & \text{if the cabling is negative},
   \end{cases}
\]
and
\[
\ve_n=
   \begin{cases}
   -1 & \text{if the cabling is negative and $q>0$}  \\
   1 & \text{otherwise}.
   \end{cases}
\]
The sign $\ve$ is chosen such that
$I(H,\gamma)>0$, where $H$ is the fibers of the Seifert fibration on $\bd(D^2\times S^1)_n$
and $I(H,\gamma)$ is the algebraic intersection number of $H$ and $\gamma$ on $\bd N$.
This is checked as follows:
$H=\ve\,p\,\fM_n+\ve\,q\,\fL_n=\ve\,q\,\fM+\ve\,p\,\fL$ on $\bd N$
and $I(H,\gamma)=\ve\,(qv-pu)$.
If the cabling is positive then we have $qv-pu>0$.
If it is negative then $qv-pu<0$. In either case, we have $I(H,\gamma)>0$.
This inequality means that $H$ intersects $F$ positively transversely,
see Figure~\ref{fig13}. 
The sign $\ve_n$ is needed since the working orientation 
of $S_n$ changes depending on the mutual positions of $0$, $q/p$ and $u/v$,
where $0$ is the slope of the longitude, $q/p$ is the slope of the cabling,
and $u/v$ is the slope of the fiber surface.

Let $\mathfrak L$ be the set of longitude $\fL$ such that $u\geq 0$ and $q\ne 0$, then
there exists a longitude $\fL$ in $\mathfrak L$ such that $u$ becomes minimal among them.
We always use this meridian-longitude pair $(\fM,\fL)$ in the discussion below.
In particular, the case $q=0$ is excluded.

\begin{lemma}\label{lemma3001}
$L(\m)\cap N$ intersects the fiber surface of $L'(\m')$ transversely.
\end{lemma}

\begin{proof}
The multilink $L(\m)\cap N$ is parallel to the fibers of the Seifert fibration of
$L_{p,q}(\m_{p,q})$ in $S^3$, denoted by $H$. So, it is enough to show that
$H$ is transverse to the fiber surface of $L'(\m')$.
By \cite[Theorem~4.2]{en},
the fibration of $L(\m)$ is decomposed into two fibered multilinks $L'(\m')$ and $L_{p,q}(\m_{p,q})$ 
by the splice decomposition, each of whose fibration is induced from that of $L(\m)$.
So, $H$ is transverse to the fiber surface of $L'(\m')$ if and only if 
$H$ is transverse to the fiber surface of $L_{p,q}(\m_{p,q})$.
We always have this transversality since $L_{p,q}(\m_{p,q})$ is fibered.
\end{proof}

\begin{lemma}\label{lemma3000}
For each $i=1,\ldots,n-1$, $m_i>0$ if and only if the cabling is positive.
\end{lemma}

\begin{proof}
Recall that the orientation of $m_iS_i$ is consistent with that of $\fL$.
If the cabling is positive then the working orientation of $S_i$ is consistent with that
of $\fL$. Hence $m_i>0$. If it is negative then, since we change the orientation of the fibers 
of the Seifert fibration by multiplying $\ve$, 
the working orientation becomes opposite to that of $\fL$. Hence $m_i<0$.
\end{proof}

\subsection{Proof of Theorem~\ref{thm04}}

\begin{lemma}\label{thm3}
Let $L(\m)$ be a fibered multilink in an oriented, closed, smooth $3$-manifold $M$ 
with a cabling in a solid torus $N$. Then there exists a positive contact form $\alpha$ on $M$
with the following properties:
\begin{itemize}
\item[(1)] $L(\m)$ is compatible with the contact structure $\xi=\ker\alpha$.
\item[(2)] On a neighborhood of $\bd N$, $\alpha$ is given as $\alpha=h_2(r)d\mu-h_1(r)d\lambda$
such that $u/v-h_1(1)/h_2(1)>0$ is sufficiently small,
where $(r,\mu,\lambda)$ are coordinates of $N=D^2\times S^1$ chosen such that
$(r,\mu)$ are the polar coordinates of $D^2$ of radius $1$ and
$(\mu,\lambda)$ are the coordinates of $\bd N$ with respect to
the meridian-longitude pair $(\fM,\fL)$, 
and $h_1$ and $h_2$ are real-valued smooth functions with parameter $r\in [0,1]$.
\item[(3)] $\alpha$ on $N$ is the restriction of the contact form compatible with
the Seifert multilink $L_{p,q}(\m_{p,q})$ to $S^3\setminus\text{{\rm int}} N(S_n)$.
\end{itemize}
\end{lemma}

\begin{proof}
Let $L'(\m')$ be the multilink in $M$ before the cabling and let $\alpha'$ be
a contact form obtained in Proposition~\ref{lemmatw}, whose kernel is compatible with $L'(\m')$.
On a neighborhood of $\bd N$, $\alpha'$ is given as
\[
   \alpha'=Rvd\mu+\left(\frac{1}{r}-Ru\right)d\lambda,
\]
as in equation~\eqref{eqsec3}. Hence
\[
\frac{u}{v}-\frac{h_1(1)}{h_2(1)}=\frac{u}{v}-\frac{-(1-Ru)}{Rv}=\frac{1}{Rv}>0
\]
can be sufficiently small since we can choose $R>0$ sufficiently large.

Next we make a contact form compatible with $L(\m)$ from $\alpha'$ by replacing the form on $N$
suitably. Let $\alpha_{p,q}$ be a positive contact form on $S^3$ 
whose kernel is compatible with the fibered Seifert multilink $L_{p,q}(\m_{p,q})$ of the cabling.
Such a contact form is given explicitly in Proposition~\ref{thm1} and Proposition~\ref{thm2}.
Let $(r_n,\mu_n,\lambda_n)$ be the coordinates on $(D^2\times S^1)_n$,
then in a small neighborhood of $\bd N$, the gluing map of the cabling is given as
$(r,\mu,\lambda)=(2-r_n,\lambda_n,\mu_n)$. Hence, on this neighborhood, we have
\[
\alpha=h_2(r)d\mu-h_1(r)d\lambda=-h_1(2-r_n)d\mu_n+h_2(2-r_n)d\lambda_n.
\]

First consider the case where the cabling in $N$ is positive.
In this case, we have $H=\ve q\,\fM_n+\ve p\,\fL_n=q\,\fM+p\,\fL$
since $\ve=1$, $q>0$, $u\geq 0$, $v>0$ and $qv-pu>0$.
By choosing $R>0$ sufficiently large, we can assume that
$H$, $\gamma$, $\alpha'$ and $\alpha_{p,q}$ are as shown in Figure~\ref{fig14}. 
Remark that the contact forms $\alpha'$ and $\alpha_{p,q}$ 
in the figures are given with the coordinates $(r_n, \mu_n, \lambda_n)$,
so the $x$-axis represents $-h_2(2-r_n)$ and the $y$-axis does $-h_1(2-r_n)$.
By multiplying a positive constant to $\alpha_{p,q}$ if necessary,
we can connect the two contact forms $\alpha'$ and $\alpha_{p,q}$ smoothly
with keeping the positive transversality
of the Reeb vector field and the interiors of the fiber surfaces.

\begin{figure}[htbp]
   \centerline{\input{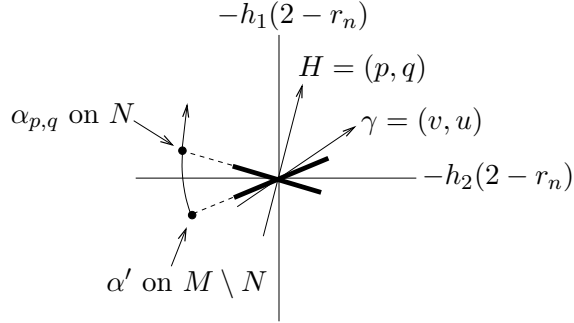}}
   \caption{Connect $\alpha'$ and $\alpha_{p,q}$ smoothly
(case of positive cabling).\label{fig14}}
\end{figure}

Next we consider the case where the cabling is negative.
Recall that the contact form constructed
according to Lemma~\ref{lemma2} and Proposition~\ref{thm2}
depends on the choice of $b_1,\ldots,b_k$.
By Lemma~\ref{lemma3000} we have $m_i<0$ for $i=1,\ldots,n-1$.
We now choose for instance $m_1S_1$ as
the negative component with index $i_0$ specified in Lemma~\ref{lemma2}.
In this setting, we re-choose these $b_i$'s such that $b_n/a_n\leq 0$, and 
then choose $R_n$ in Lemma~\ref{lemma2}~(2) sufficiently large so that
the line representing $\ker\alpha_{p,q}$ is sufficiently close to $H$ on the $xy$-plane.

If $q<0$ then we have $H=\ve q\,\fM_n+\ve p\,\fL_n=q\,\fM+p\,\fL$ since $\ve=1$,
$u\geq 0$, $v>0$ and $qv-pu>0$.
By choosing $R>0$ sufficiently large, we can assume that
$H$, $\gamma$, $\alpha'$ and $\alpha_{p,q}$ are as shown on the left in Figure~\ref{fig14b}.
If $q>0$ then $H=\ve q\,\fM_n+\ve p\,\fL_n=-q\,\fM-p\,\fL$ since $\ve=-1$.
Thus they are in the positions as shown on the right in Figure~\ref{fig14b}.
In either case,
by multiplying a positive constant to $\alpha_{p,q}$ if necessary,
we can connect the contact forms $\alpha'$ and $\alpha_{p,q}$ smoothly as shown in
the figures. Thus we obtain the contact form required.
\end{proof}

\begin{figure}[htbp]
   \centerline{\input{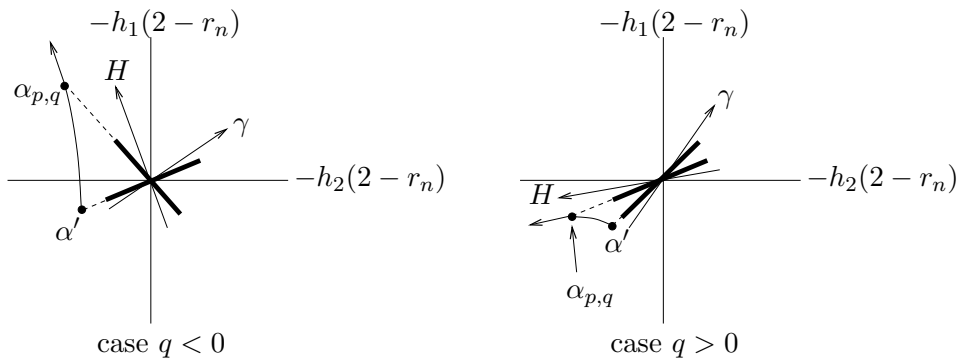}}
   \caption{Connect $\alpha'$ and $\alpha_{p,q}$ smoothly
(case of negative cabling).\label{fig14b}}
\end{figure}

Now we prove Theorem~\ref{thm04}. We first recall the statement.
\vspace{3mm}

\noindent
{\bf Theorem~\ref{thm04}.\;}
Let $L(\m)$ be a fibered multilink in an oriented, closed, smooth $3$-manifold $M$ 
with cabling in a solid torus $N$ in $M$ and $L'(\m')$ be the fibered multilink 
obtained from $L(\m)$ by retracting $N$ into its core curve.
Let $\xi$ and $\xi'$ denote the contact structures on $M$ compatible with 
$L(\m)$ and $L'(\m')$ respectively.
\begin{itemize}
\item[(1)] If $\xi'$ is tight and the cabling is positive, then $\xi$ is tight.
\item[(2)] If $\xi'$ is tight, the cabling is negative and $L(\m)\cap N$ has at least two components,
then $\xi$ is overtwisted.
\item[(3)] If $\xi'$ is tight, the cabling is negative, $L(\m)\cap N$ is connected, $p\geq 2$ and $q\leq -2$,
then $\xi$ is overtwisted.
\item[(4)] If $\xi'$ is overtwisted then $\xi$ is also overtwisted.
\end{itemize}

\begin{proof}
We use the contact structure constructed in Lemma~\ref{thm3}.
If $\xi'$ is in case~(1) in the assertion then there exists a one-parameter family
which connects $\xi$ and $\xi'$. Hence $\xi$ and $\xi'$ are contactomorphic by Gray's theorem~\cite{gray}. 
Suppose that $\xi'$ is in case~(2). In this case, each $m_iS_i$ for $i=1,\ldots,n-1$ is 
a negative component of $L_{p,q}(\m_{p,q})$ by Lemma~\ref{lemma3000}. Thus,
Proposition~\ref{thm2} and Lemma~\ref{thm3} ensure that 
there exists a negative component which contains an overtwisted disk.
Suppose $\xi'$ is in case~(3).
We will use Lemma~\ref{lemma3} to detect an overtwisted disk.
We assign the index $i_0$ to the link component $S_1$ and the index $i_1$ to the singular fiber 
of the Seifert fibration other than $S_n$.
From Figure~\ref{fig14}, we can make sure that
the proof of Lemma~\ref{thm3} works even if the point representing $\alpha'$ is sufficiently
close to the horizontal axis.
This means that we can choose $R_n$ to be any value in $(-b_n/a_n, \infty)$.
This is important since, in the proof of Lemma~\ref{lemma3}, 
$R_n$ is some value with $-b_n/a_n<R_n$ and we do not know
at which value the overtwisted disk is detected. Since $a_{i_0}=1$, we have $1>1/|p|+1/|q|$.
So, we can detect an overtwisted disk between 
$(D^2\times S^1)_{i_0}$ and $(D^2\times S^1)_{i_1}$ by Lemma~\ref{lemma3}, 
which is outside of $(D^2\times S^1)_n$.
In case~(4), let $D$ denote an overtwisted disk in $(M,\xi')$.
Since we can choose $N$ sufficiently small such that $\bd D\cap N=\emptyset$,
the overtwisted disk still remains in $(M,\xi)$ after the cabling.
\end{proof}

\begin{rem}
(1)\;\,If $p=1$ then $L(\m)$ is ambient isotopic to $L'(\m')$. 
Suppose $p\geq 2$.
We have chosen $(\fM,\fL)$ such that $u\geq 0$ is minimal among $\mathcal L$.
If the cabling is negative and $q\geq 2$ then we can change $\fL\mapsto \fL-(q-1)\,\fM$ such that 
the cabling is negative and $q=1$. Hence
this case is excluded since $u$ is not minimal in $\mathcal L$.
Now, the remaining case becomes when $\xi'$ is tight, the cabling is negative, 
$L(\m)\cap N$ is connected, $p\geq 2$ and $q\in\{-1, 1\}$.\\
(2)\;\,
We have excluded the case $q=0$. This is because we only gave explicit constructions of 
contact forms when $A\ne 0$.
Actually, it is not difficult to give a contact form with the same property explicitly when $A=0$, i.e., $q=0$.
If we include the case $q=0$ in the above argument,
the remaining case becomes when $\xi'$ is tight, the cabling is negative, 
$L(\m)\cap N$ is connected, $p\geq 2$ and $q\in\{-1, 0\}$.
\end{rem}

\subsection{Cabling along fibered knots}

Let $L'$ be a fibered knot in $M$ and $N(L')$ its small, compact, tubular neighborhood
with the meridian-longitude pair $(\fM',\fL')$ determined by the fiber surface, namely
$\fM'$ is the boundary of a meridional disk
and $\fL'$ is the oriented boundary of a fiber surface of $L'$.

\begin{cor}\label{cor7}
Let $L'$ be a fibered knot in an oriented, closed, smooth $3$-manifold $M$ 
and $L$ be the link obtained from $L'$ by cabling a $(p,q)$-torus link with respect to
$(\fM',\fL')$, i.e., the cabling with slope $q\,\fM'+p\,\fL'$.
Let $\xi$ and $\xi'$ denote the contact structure on $M$ compatible with 
$L$ and $L'$ respectively. 
\begin{itemize}
\item[(1)] If $\xi'$ is tight and $q>0$ then $\xi$ is tight.
\item[(2)] If $\xi'$ is tight, $q<0$ and $\gcd(p,|q|)\geq 2$ then $\xi$ is overtwisted.
\item[(3)] If $\xi'$ is tight, $p\geq 2$ and $q\leq -2$ then $\xi$ is overtwisted.
\item[(4)] If $\xi'$ is overtwisted then $\xi$ is also overtwisted.
\end{itemize}
\end{cor}

\begin{proof}
Let $L'(\m')$ be the fibered multilink obtained from $L$ by
retracting the solid torus $N(L')$ of the cabling to its core curve. 
Since $L'$ is a knot, the framing of the fiber surfaces of 
$L'(\m')$ is given by the boundary of a fiber surface of $L'$. This means that
$\gamma=\fL$, i.e., $(u,v)=(0,1)$.
Hence the cabling is positive in the sense in Definition~\ref{dfn04}
if and only if $q>0$. Note that the case $q=0$ is excluded by Lemma~\ref{lemma3001}.
Then, the assertion is just a restatement of Theorem~\ref{thm04} in this special case.
\end{proof}

\begin{rem}
It is known in~\cite{behm} that  
in the remaining case, i.e., the case where $\xi'$ is tight, $p\geq 2$ and $q=-1$,
the contact structure $\xi$ is tight if and only if $M=S^3$ and $L$ is a trivial knot
(cf.~\cite{hedden2} for the case where $L'$ is a fibered knot in $S^3$).
\end{rem}


\begin{thebibliography}{99} 
\bibitem{bi}S.~Baader and M.~Ishikawa,
   {\it Legendrian graphs and quasipositive diagrams}, 
   Ann. Fac. Sci. Toulouse Math. (6) {\bf 18} (2009), 285--305.
\bibitem{behm}K.~Baker, J.~Etnyre and J.~van Horn-Morris,
   {\it Cablings, contact structures and mapping class monoids},
   preprint, arXiv:1005.1978 [math.SG].
\bibitem{bm}G.~Burde and K.~Murasugi,
   {\it Links and Seifert fiber spaces},
   Duke Math. J. {\bf 37} (1970), 89--93.
\bibitem{ch}V.~Colin and K.~Honda,
   {\it Reeb vector fields and open book decompositions}
   preprint, available at arXiv:0809.5088v1 [math.GT].
\bibitem{en}D.~Eisenbud and W.~Neumann,
   {\it Three-dimensional link theory and invariants of plane curve singularities},
   Ann. of Math. Stud. 110, Princeton Univ. Press, Princeton, NJ, 1985.
\bibitem{eliashberg:89}
   Y.~Eliashberg,
   {\it Classification of overtwisted contact structures on 3-manifolds},
   Invent. Math. {\bf 98} (1989), 623--637.
\bibitem{eliashberg:90}
   Y. Eliashberg,
   {\it Filling by holomorphic discs and its applications},
   Geometry of low-dimensional manifolds, 2 (Durham, 1989), pp. 45--67,
   London Math. Soc. Lecture Note Ser., 151, Cambridge Univ. Press, Cambridge, 1990.
\bibitem{eliashberg}
   Y. Eliashberg,
   {\it Contact 3-manifolds twenty years since J. Martinet's work},
   Ann. Inst. Fourier (Grenoble) {\bf 42} (1992), 165--192.
\bibitem{eo}T.~Etg\"{u} and B.~Ozbagci,
   {\it Explicit horizontal open books on some plumbings},
   Internat. J. Math. {\bf 17} (2006), 1013--1031.
\bibitem{etnyre}J. B.~Etnyre,
   {\it Lectures on open book decompositions and contact structures},
    Floer homology, gauge theory, and low-dimensional topology, pp. 103--141, 
    Clay Math. Proc., 5, Amer. Math. Soc., Providence, RI, 2006.
\bibitem{geiges}H.~Geiges,
   {\it An Introduction to Contact Topology},
   Cambridge Stud. in Adv. Math. 109, Cambridge Univ. Press, Cambridge, 2008.
\bibitem{giroux}E.~Giroux,
   {\it G\'eom\'etrie de contact:
   de la dimension trois vers dimensions sup\'erieures},
   Proceedings of the International Congress of Mathematicians,
   Vol II (Beijing, 2002), pp.405--414,
   Higher Ed. Press, Beijing, 2002.
\bibitem{gromov}M.~Gromov,
   {\it Pseudoholomorphic curves in symplectic manifolds}, 
   Invent. Math. {\bf 82} (1985), 307--347.
\bibitem{gray}J.~Gray,
   {\it Some global properties of contact structures},
   Ann. of Math. {\bf 69} (1959), 421--450.
\bibitem{hedden} M.~Hedden,
   {\it Notions of positivity and the Ozsv\'ath-Szab\'o concordance invariant},
   J. Knot Theory Ramifications {\bf 19} (2010), 617--629.
\bibitem{hedden2} M.~Hedden,
   {\it Some remarks on cabling, contact structures, and complex curves},
   Proc. G\"{o}kova Geom. Topol. Conference, 2007, pp. 49--59,
   G\"{o}kova Geom./Topol. Conference, G\"{o}kova, 2008.
\bibitem{ishikawa3} M.~Ishikawa,
   {\it Compatible contact structures of fibered positively-twisted graph multilinks in the 3-sphere},
   preprint, available at arXiv:1006.4414 [math.GT].
\bibitem{kt} Y.~Kamishima and T.~Tsuboi,
   {\it CR-structures on Seifert manifolds},
   Invent. Math. {\bf 104} (1991), 149--163.
\bibitem{lm}P.~Lisca and G.~Mati\'{c},
   {\it Transverse contact structures on Seifert $3$-manifolds},
   Algebr. Geom. Topol. {\bf 4} (2004), 1125--1144.
\bibitem{mw}J. D.~McCarthy and J. G.~Wolfson,
   {\it Symplectic gluing along hypersurfaces and resolution of isolated orbifold singularities},
   Invent. Math. {\bf 119} (1995), 129--154.
\bibitem{ozbagci}B.~Ozbagci,
   {\it Explicit horizontal open books on some Seifert fibered $3$-manifolds},
   Topology Appli. {\bf 154} (2007), 908--916.
\bibitem{os}B.~Ozbagci and A. I.~Stipsicz,
   {\it Surgery on Contact $3$-manifolds and Stein Surfaces},
   Bolyai Soc. Math. Stud. 13, Springer-Verlag, Berlin, 2004.
\bibitem{rudolph1}L.~Rudolph,
   {\it A characterization of quasipositive Seifert surfaces 
   (Constructions of quasipositive knots and links, III)},
   Topology {\bf 31} (1992), 231--237.
\bibitem{rudolph2}L.~Rudolph,
   {\it Quasipositive annuli (Constructions of quasipositive knots and links, IV)},
    J. Knot Theory Ramifications {\bf 1} (1992), 451--466.
\bibitem{rudolph3}L.~Rudolph,
    {\it Quasipositivity as an obstruction to sliceness},
    Bull. Amer. Math. Soc. (N.S.) {\bf 29} (1993), 51--59. 
\bibitem{rudolph4}L.~Rudolph,
   {\it Quasipositive plumbing (constructions of quasipositive knots and links, V)},
    Proc. Amer. Math. Soc. {\bf 126} (1998), 257--267.
\bibitem{rudolph5}L.~Rudolph,
    {\it Positive links are strongly quasipositive},
    Proceedings of the Kirbyfest (Berkeley, CA, 1998), pp. 555--562,
    Geom. Topol. Monogr., 2, Geom. Topol. Publ., Coventry, 1999.
\bibitem{rudolph6}L.~Rudolph,
    {\it Quasipositive pretzels},
    Topology Appl. {\bf 115} (2001), 115--123.
\bibitem{st}A.~Sato and T.~Tsuboi,
   {\it Contact structures on closed manifolds fibered by the circles},
   Mem. Inst. Sci. Tech. Meiji Univ. {\bf 33} (1994), 41--46.
\bibitem{tanaka}T.~Tanaka,
   {\it Maximal Bennequin numbers and Kauffman polynomials of positive links},
   Proc. Amer. Math. Soc. {\bf 127} (1999), 3427--3432.
\bibitem{tw}W. P.~Thurston and H.~Winkelnkemper,
   {\it On the existence of contact forms},
   Proc. Amer. Math. Soc. {\bf 52} (1975), 345--347.
\bibitem{torisu}I.~Torisu,
   {\it Convex contact structures and fibered links in 3-manifolds},
   Internat. Math. Res. Notices {\bf 9} (2000), 441--454.
\end{thebibliography}
\end{document}